\documentstyle[epsfig,11pt]{article}
\input epsf
\onecolumn 
\begin{document}



\newcommand {\RR}{\,\hbox{\rm R}\!\!\!\!\!\hbox{\rm I}\,\,\,}

\newcommand{\square}{\hbox{${\vcenter{\hrule height.4pt
  \hbox{\vrule width.4pt height6pt \kern6pt
     \vrule width.4pt}
  \hrule height.4pt}}$}}
\newcommand {\ProofEnd} {\hfill \nobreak \square \medbreak}

\newcommand {\Hess} [3] {Hess\,{#1} \left( {#2}, {#3} \right)}
\newcommand {\Pind} {\hspace{.51cm}}  
\newcommand {\Dminr}{\bigg(\frac{Vol(b^{-1}(r_1)}{w_n 2^n}\bigg)^{(1/(n-1))}}
\newcommand {\Dmin}{\bigg(\frac{V_\infty/2}{w-n 2^n}\bigg)^{(1/(n-1))}}
\newcommand {\Om} {\Omega}
\newcommand {\vare}{\varepsilon}
\newcommand {\grad}{\bigtriangledown}
\newcommand {\Ftil}{\tilde{{\cal{F}}}}  
\newcommand {\F}{{\cal{F}}}
\newcommand {\Gtil}{\tilde{{\cal{G}}}}  
\newcommand {\G}{{\cal{G}}}



\newtheorem {theorem} {Theorem}
\newtheorem {proposition} [theorem] {Proposition}
\newtheorem {lemma} [theorem] {Lemma}
\newtheorem {definition} [theorem] {Definition}
\newtheorem {corollary} [theorem] {Corollary}
\newtheorem {remark} [theorem] {Remark}
\newtheorem {note} [theorem] {Note}
\newtheorem {example} [theorem] {Example} 



%
%
%

\begin{titlepage} 
\begin{center}
{\huge\bf The Almost Rigidity of Manifolds with\\  
\vspace*{2ex} Lower Bounds on Ricci Curvature\\
\vspace*{2ex} and Minimal Volume Growth}

\vspace*{7ex} 
{\Large Christina Sormani}$ ^\dagger$\footnote{
This material is based upon work supported under a National Science
Foundation Graduate Fellowship.\\
        \hspace*{1em}$ ^\dagger$The author can be reached by e-mail at
       {\tt sormani@math.jhu.edu}\\ 
       {\it Department of Mathematics\\
        Johns Hopkins University\\
        3400 North Charles Street\\
        Baltimore MD 21218}\\} 

\vspace*{7ex} 
{\large To appear in Communications in Analysis and Geometry}
\date{}
\end{center}
\end{titlepage} 
\pagestyle{plain}  
\pagenumbering{arabic}



\newpage

\Pind
Twenty years ago, Calabi and Yau each proved that 
a complete noncompact Riemannian manifold with nonnegative Ricci curvature 
must have at least linear volume growth [Yau].  This was proven by studying the
Busemann function, $b=b_\gamma$, associated with a ray, $\gamma$,
$$
b(x)=\lim_{R\to\infty} ( R-d(x,\gamma(R)).
$$
In [So2], the author proved that if such a manifold has linear volume growth
then its Busemann functions are proper.  The simplest examples of manifolds 
with linear volume growth  are the metric product manifolds, $X\times \RR$,
whose  cross sections, $X\times \{ r \}$, are 
level sets of the Busemann functions. 

In this paper we prove that a complete
noncompact manifold with nonnegative Ricci curvature and linear volume
growth must be close to being such a metric product manifold asymptotically
[Theorem~\ref{GHCross}].  That is, as $r \to \infty$, the set
$b^{-1}([r, r+L])$ becomes close to $b^{-1}(r) \times [r, r+L]$ in
the Gromov-Hausdorff topology where the closeness depends linearly 
on $diam(b^{-1}(r))$.  See Section 2.  The proof 
involves a  careful analysis of the Busemann function using the 
recently-developed 
Cheeger-Colding Almost Rigidity Theory [ChCo].  We also use this method to
prove the following theorem.

\begin{theorem} \label{SublinDiam2}
If $M^n$ is a manifold with nonnegative Ricci curvature and linear
volume growth, then it has sublinear diameter growth,
$$
\lim_{R\to\infty} \frac {diam(b^{-1}(R))}{R}=0.
$$
\end{theorem}

Examples of manifolds satisfying the hypothesis of this theorem
for which $diam(b^{-1}(R))$ grows logarithmicly appeared in [So2].
Applications of this theorem to the analysis
of harmonic functions on manifolds with nonnegative Ricci curvature
and linear volume growth will appear in [So3].

In addition to studying manifolds with nonnegative Ricci curvature, 
we study manifolds with a quadratically decaying lower bound
on Ricci curvature,
\begin{equation} \label{RicciBound}
Ric(x) \ge \frac {(n-1)(\frac 1 4 - v^2)} {b(x)^2}
\textrm{ for all } x \in b^{-1}([r_0,\infty)),
\end{equation}
where $v$ can take any value in 
$[0,\frac {n+1}{2(n-1)})$.
This bound is implied by the traditional quadratically decaying lower bound
on Ricci curvature defined using the distance function from a base point
[ChGrTa].

Such manifolds also have a lower bound on their volume growth,
\begin{equation} 
\liminf_{R\to\infty}\frac {Vol(B_{x_0}(R))}{R^p} \ge C >0,
\end{equation}
where $p=(n-1)(1/2-v)+1>0$ [ChGrTa].  If one is given the (n-1) Hausdorff
volume of a compact set, $S \subset b^{-1}(r_1) \subset b^{-1}((r_0,\infty))$,
then there is a precise lower bound, $C=C_S$, depending upon that volume
[So2, see Corollary~\ref{TheoremI}].  In section 1.1, we review this 
and a relative volume comparison theorem for Busemann functions that is 
needed for this paper.

\begin{definition} \label{DefnMin}\label{MinVolGrowth}  
{\em We say that a manifold, $M^n$, has {\em minimal volume growth} if 
it has a quadratically decaying lower Ricci curvature bound as in
(\ref{RicciBound}) and }
$$
\limsup_{R\to\infty}\frac {Vol(B_{x_0}(R))}{R^p} = V_0 < \infty,
$$
{\em where $p=(\frac 1 2 -v)(n-1)+1$.}
\end{definition}

Note that a manifold with nonnegative Ricci curvature, $v=1/2$, has minimal
volume growth if it has linear volume growth.

In contrast, we will say that  a manifold 
satisfying (\ref{RicciBound}) has 
{\em strongly minimal volume growth} with respect to the compact set, 
$S\subset b^{-1}(r_1)$, if
\begin{equation} \label{Mstrongmin}
\limsup_{R\to\infty}
\frac {Vol(B_{x_0}(R)))}{R^p} = C_S
\end{equation}
where $p=(\frac 1 2 -v)(n-1)+1$ [see also Defn~\ref{StronglyMin}].

In Section 1.2 we state and prove a {\em rigidity theorem} 
for the end of a manifold
with strongly minimal volume growth [Theorem~\ref{VolMinEndWarp}].  First we
show that the given level set,  $b^{-1}(r_1)$, is  
a compact smoothly embedded submanifold.
Then we prove the end of the manifold is a warped product,
$$
b^{-1}([r_1,\infty))=
b^{-1}(r_1) \, \times_{(b/r_1)^{(1/2 -v)}} \, [r_1, \infty).
$$ 
That is, it has a metric of the form
\begin{equation} \label{WarpProdMet}
db^2\, +\, \Big(\frac {b} {r_1}\Big)^{(\frac 1 2 -v)2}g_0
\end{equation}
where $g_0$ is the induced metric on $b^{-1}(r_1)$.

In particular, a manifold with $Ricci \ge 0$, $v=1/2$ and strongly
minimal volume growth has only one end and that end is an isometric
product [Cor~\ref{RigidIsom}].

In [So2], there are examples which demonstrate that the 
{\em strongly minimal volume growth} condition in Theorem~\ref{VolMinEndWarp}
is necessary.  
These examples have nonnegative Ricci curvature and linear volume growth,
$$
\lim_{r\to\infty}\frac {Vol(B_{x_0}(r))}{r} = V_0 < \infty,
$$
but their ends are not isometric product manifolds.

In the second section of this paper, we prove the 
corresponding {\em almost rigidity theorem},
Theorem~\ref{GHWarp}.  Here we add the assumption that $Ricci \ge 0$
everywhere on the manifold along with the quadratically decaying lower
bound on Ricci curvature, (\ref{RicciBound}), with $v$ now in $[0, 1/2]$
and minimal volume growth.  These conditions imply that the Busemann
function is proper [So2, Theorem 19].  

Theorem~\ref{GHWarp} roughly states that such a manifold
is asymptotically close in the Gromov-Hausdorff sense 
to a warped product manifold with a metric as in (\ref{WarpProdMet}).  
More precisely, we show that for $r$ sufficiently large, the compact
region $b^{-1}([r, r+L])$ is close to a warped product, 
$X_r \times_f [r, r+L]$ with the warping function, $f(s)=(s/r)^{(1/2-v)}$, 
where $X_r$ is
a length space close to $b^{-1}(r)$,
[Theorem~\ref{GHWarp}].  This closeness depends linearly
on the diameter $b^{-1}(r)$.  

It is important to note that a manifold with minimal volume growth and 
with globally nonnegative and quadratically decaying 
Ricci curvature is not asymptotically close 
to a unique warped product manifold.  In particular a manifold with 
nonnegative Ricci curvature and linear volume growth has regions, 
$b^{-1}((r, r+L))$,
which approach isometric product manifolds, $X_r \times (r, r+L)$,
but $X_r$ may change slowly as $r$ approaches infinity.  Examples where
$X_r$ alternates between two different Riemannian manifolds appear in [So2].

In order to precisely state Theorem~\ref{GHWarp}, we must define
{\em localized distance}, the metric which we will be using
on the regions, $b^{-1}([r, r+L])$.  This definition and a definition of
the Gromov-Hausdorff distance appear in Section 2.1.  The statement
of Theorem~\ref{GHWarp} and the related Theorem~\ref{GHCross}, along 
with a discussion
of their consequences appears in Section 2.2.  The proof appears in 
Section 2.4 after a Section 2.3, which relates Cheeger and Colding's 
work on almost rigidity and maximal volume to our condition of minimal 
volume growth.

Lastly, in Section 3 we focus on manifolds with globally nonnegative
Ricci curvature, $v=1/2$, and linear volume growth.
The same lemmas used to prove our almost rigidity theorem in Sections 2.1 
and 2.3 are directly applied to prove two results on the diameter 
growth of such manifolds.  In
Theorem~\ref{SublinDiam}, we prove that the {\em localized diameter} of 
the level sets of the Busemann function grows sublinearly. 
This theorem is stated in Section 3 after defining {\em localized diameter}
in Definition~\ref{DefnDiam}.  This theorem is then employed to prove the 
Theorem~\ref{SublinDiam2} stated above, in which the diameter of 
the level sets are measured in the ambient manifold, $M^n$.

The author would like to thank Professor Cheeger for suggesting 
an analysis of manifolds with minimal volume growth and for numerous
enlightening conversations.  She would also like to thank Professor Colding
for helpful discussions on his work in almost rigidity theory. 
Finally, she is very grateful to the Courant Institute
of Mathematical Sciences for its generous support in her years
as a graduate student.
 
Background material can be found in [ChEb], [BiCr] and [Ch].

\section{Strongly Minimal Volume Growth and Rigidity}
%
%

\subsection{Background on Volume Growth}

We now review some definitions and theorems regarding special sets in 
noncompact
manifolds with quadratically decaying lower Ricci curvature bounds.  
See [So2].  This background will be used to study both the 
rigidity and almost rigidity of such manifolds later in this paper.  

All geodesics and rays are parametized by arclength.

\begin{definition} \label{DefnRay} {\em
Let $\gamma$ be a ray, $b_\gamma(x)$ be its associated Busemann function and
$x\in M^n$.   A ray, $\gamma_x:[b(x),\infty) \mapsto M$, emanating from $x$ 
is called a {\em Busemann ray} associated with $\gamma$
if it is the limit of a sequence of minimal geodesics $\sigma_i$ from
$x$ to $\gamma(R_i)$ in the following sense,
$$
\gamma'_x(b(x))= \lim_{R_i\to\infty} \sigma_i'(0).
$$
Note that $\gamma_x(b(x))=x$. Note also that $\gamma_x$ need not be unique.}
\end{definition}  

\begin {definition}  \label{DefnOmega}
Given a compact set, $K$, contained in $M^n$ and a ray, $\gamma$,
let
$$
\Omega(K)=
\{x: \, \exists z \in K \, \exists t \ge b(z) \textrm{ s.t. } x=\gamma_z(t)\}
$$
Let
$$
\Omega_{s_1,s_2}(K)=\Omega(K) \cap b^{-1}([s_1,s_2]).
$$
\end{definition}

In [So2, Cor 15], the author proved the following volume comparison theorem.

\begin{theorem} [So2, Thm 5, Cor 17]  \label{Monotonicity}
Let 
$
Ric \ge  {(n-1)(\frac 1 4 - v^2)}/ {b(x)^2}
$
whenever $x\in b^{-1}([r_0,\infty))$
where $v\in [0, \frac{n+1}{2(n-1)})$.
Let $r_1 \ge r_0$ and let $K$ be a compact set contained in 
$b^{-1}((-\infty,r_1])$.

Then there exists a nondecreasing function, $V(r)$, such that
\begin{equation}
V(r)=\frac {Vol_{n-1}(\Omega(K)\cap b^{-1}(r))}{r^{(\frac 1 2 -v)(n-1)}}
\end{equation}
almost everywhere in $[r_1,\infty)$.

In particular, for almost every $s_2 > s_1 \ge r_1$ we have
\begin{equation}
\frac {V(s_1)}{p}\le 
\frac {Vol(\Omega(K)\cap b^{-1}(s_1,s_2))}{s_2^{p}-s_1^p} 
\le
\frac {V(s_2))}{p},
\end{equation}
where $p=(1/2-v)(n-1)+1$.
\end{theorem}

This theorem was proven using the following series of comparison manifolds 
which were also employed in [ChGrTa].  

\begin{definition} \label{CompMan}
The comparison manifold, $M^n_{R,\vare}$, is a warped product manifold
diffeomorphic to $\RR^n$, with the metric 
$dt^2 + J_{R,\vare}(t)^2 g_0$ where
\begin{equation}\label{Jacobi2}
J_{R,\vare}(t)=\frac {R+\vare} {2v}\left(
            -\left(\frac{R-t+\vare}{R+\vare}\right)^{\frac 1 2 +v}
              +\left(\frac{R-t+\vare}{R+\vare}\right)^{\frac 1 2 -v}\right)
\end{equation}
and $g_0$ is the standard metric on the sphere.
Here $\vare:=0$ when $v \ge 1/2$. 
\end{definition}

These manifolds satisfy the following Ricci curvature
bound,
\begin{equation} \label{RicciComp1}
Ric_y \ge \frac {(n-1)(\frac 1 4 -v^2)}{(R-d(y,0)+\vare)^2}.
\end{equation}

Note also that
\begin{equation} \label{JacobiLimit}
\lim_{R\to\infty}\frac{\int_{R-r_4}^{R-r_3}J_{R.\vare}(t)^{n-1}  \,dt }
{\int_{R-r_3}^{R-r_1} J_{R.\vare}(t)^{n-1} \,dt }
=\frac {(r_4)^p-(r_3)^p}{(r_3)^p-(r_1)^p}.
\end{equation}

These comparison manifolds, (\ref{JacobiLimit}),
and the Relative Volume Comparison Theorem [BiCr,GrLaPa] 
will be used to prove our rigidity theorem.  

The following corollary of Theorem~\ref{Monotonicity}, defines the
precise constant, $C_S$, mentioned in the introduction and allows us
to define {\em strongly minimal volume growth} for ends.

\begin{corollary}  \label{theoremI} \label{TheoremI}
Given $M^n$ as described in (\ref{RicciBound}) and given any $R>r_1\ge r_0$
and any $R_0 >0$, let $S=B_{x_0}(R_0)\cap b^{-1}(r_1)$.
Then
\begin{equation}
Vol(B_{x_0}(R+R_0-r_0)) \geq C_S \bigg(R^p-r_0^p\bigg)
\end{equation}
where $p=(\frac 1 2 -v)(n-1)+1 >0$ and
$
C_S={Vol_{n-1}(S)}/{(p r_1^{p-1})}.
$
\end {corollary}

\begin{definition} \label{StronglyMin} 
{ \em Given $r_1 \ge r_0$, we say that a region, $b^{-1}((r_1, \infty))$, in a
manifold, $M^n$, has {\em strongly minimal volume growth} 
with respect to a given ball $B_{x_0}(R_0)$,
if it has a quadratically decaying lower Ricci curvature bound as in
(\ref{RicciBound}) and }
$$
\lim_{R\to\infty}\frac {Vol(B_{x_0}(R) \cap b^{-1}((r_1, \infty)))}{R^p} = 
C_S, 
$$
{\em 
where $S=B_{x_0}(R_0)\cap b^{-1}(r_1)$
and $p=(\frac 1 2 -v)(n-1)+1$ }.
\end {definition}

\subsection{Strongly Minimal Volume Growth
and Rigidity} \label{RigidSect}

In this section we prove the following rigidity theorem.

\begin{theorem} \label{VolMinEndWarp}
Let $M^n$  
be a complete noncompact manifold with a ray $\gamma$ and
a Busemann function $b: M \to \RR $ and a 
quadratically decaying lower Ricci curvature bound as in (\ref{RicciBound}).
Suppose $b^{-1}((r_1,\infty))$ has strongly minimal volume growth 
with respect to a given ball, $B_{x_0}(R_0)$.
Then $b^{-1}(r_1)$ is a smoothly embedded submanifold 
contained in $\bar{B}_{R_0}(\gamma(0))$ and 
$$
b^{-1}([r_1,\infty))=
b^{-1}(r_1) \, \times_{(b/r_1)^{(1/2 -v)}} \, [r_1, \infty).
$$ 
\end {theorem}

Note that we are only prescribing the manifold's properties on one
end or on a subset of that end, $b^{-1}([r_1, \infty))$,
and we only prove that that end
is rigid.  The rest of the manifold can have larger volume growth and
any kind of curvature.  

When $v=1/2$, this theorem combined with [So2, Cor 23]
and some simple calculations imply the following corollary.

\begin{corollary} \label{RigidIsom}
Let $M^n$ be a complete noncompact manifold with 
globally nonnegative Ricci curvature such that
$$
\lim_{r\to\infty}\frac {Vol(B_{x_0}(r))}{r} = 
Vol_{n-1}\bigg(B_{x_0}(R_0) \cap b^{-1}(r_1)\bigg).
$$
Then $M^n$ has only one end and that end is an isometric product
manifold, $b^{-1}(r_1) \times [r_1, \infty)$ where $b^{-1}(r_1)
\subset \bar{B}_{x_0}(R_0)$.
\end{corollary}

It is important to note that the Busemann function, in general, is
not a smooth function, although its gradient is $1$ almost everywhere.
In order to prove that the level sets are smooth
on a manifold with strongly minimal volume growth, we will show that
the Busemann function satisfies the following elliptic partial differential
equation in the weak sense:
\begin{equation}\label{eqnlapb}
\Delta b = \frac {(n-1)(\frac 1 2 -v)}{b}.
\end{equation}
This equation is satisfied by the function $b$ on a manifold with
a metric of the form
\begin{equation}\label{eqnwarpb}
db^2 + \Bigg(\frac b {r_0}\Bigg)^{(\frac 1 2 -v)2}g_0.
\end{equation}
Recall that we plan to prove that the end, $b^{-1}([r_1, \infty))$
has a metric of this form.  We can then use elliptic regularity to 
obtain smoothness.

We first show that manifolds with strongly minimal volume growth
have proper Busemann functions and 
the ratios, $Vol_{n-1}(b^{-1}([r_2,r_3]))/(r_3^p-r_2^p)$, are constant 
with respect to $r_2$ and $r_3$ [Lemmas~\ref{EndinOm} and~\ref{ConstRatio}].  
We then compare
these regions to subsets of annuli about $\gamma(R)$ using the fact
that the Busemann function's level sets are compact [Defn~\ref{InnerAnn},
Lemma~\ref{BuseInnerAnn}].  Recall that
$$
b(x)=\lim_{R\to\infty}R - \rho_R(x) \textrm{ where } \rho_R(x)=d(x,\gamma(R)).
$$ 
We then employ the Relative Volume Comparison
Theorem of [BiCr] and [GrLaPa] combined with the volume estimate
on these regions, to control $\Delta \rho_R$ in a weak
sense [Lemma~\ref{DiffA_R(t)}].  Then taking $R$ to infinity we show
that $b(x)$ satisfies the 
elliptic equation (\ref{eqnlapb}) in a weak sense [Lemma~\ref{=CompRay2}].  
Once we have proven that the 
Busemann function is smooth, we use the Bochner-Weitzenboch 
Formula to prove 
that the metric is rigid [Lemma~\ref{LapB5}].  

\begin{lemma} \label{EndinOm}
For $M^n$ as above, let $K=\bar{B}_{x_0}(R_0) \cap b^{-1}((-\infty,r_1])$.
Then
$$
b^{-1}([r_1,\infty)) \subset \Omega(K)
$$
and thus the Busemann function, $b$, is proper.
\end{lemma}

\noindent {\bf Proof:}  
By [So2, Lemma 4], $\Omega(K)\cap b^{-1}([r_1,\infty))$ is a closed set.
Thus if there exists a point in $b^{-1}([r_1,\infty))$ which is not
in $\Omega(K)$ then there exists a a ball 
$$
B_q(\delta) \subset b^{-1}([r_1, \infty)) \setminus \, \Omega(K).
$$
Let $s=b(q)+\delta/2$ and $t=b(q)-\delta$.

Let $U= b^{-1}([t,s]) \cap \bar{B}_q(\delta)$.  Thus
\begin{equation} \label{VolSlice}
Vol_{n-1}(\Omega(\bar{U})\cap b^{-1}(s))\ge
Vol_{n-1}(\bar{B}_q(\delta)\cap b^{-1}(s))>0.
\end{equation}
Applying Theorem~\ref{Monotonicity} we have,
\begin{equation} \label{undebound}
Vol(\Omega(\bar{U})\cap b^{-1}([s,r])) \ge \frac{r^p-s^p}{ps^{p-1}}  
     Vol_{n-1}(\Omega(\bar{U})\cap b^{-1}(s))
\end{equation}
Let $r_2$ be a real number such that
$$ 
B_{x_0}(r_2) \supset (\Omega(\bar{U})\cap b^{-1}(s)) \cup 
(\Omega(K)\cap b^{-1}(s)).
$$
By the definition of $\Omega(K)$ [Defn~\ref{DefnOmega}], for all $p \in
\Omega(K)\cap b^{-1}([s,t])$ there exists $z\in K$ such that 
$\gamma_z(b(p))=p$.  So, 
$$
d(\,p\, , \,\Omega(K)\cap b^{-1}(s)\,) \le b(p)-s \le t-s.
$$
Thus,
$$
B_{x_0}(r+r_2-s) \cap b^{-1}([r_1,\infty))
\supset \big(\Omega(K)\cap b^{-1}(s,r))\big)
   \cup \big(\Omega(\bar{U})\cap b^{-1}(s,r))\big)
$$
where the union is disjoint.
So by Theorem~\ref{Monotonicity} we have,
\begin{eqnarray}
\lefteqn{
\lim_{r\to\infty}\frac{Vol(B_{x_0}(r+r_2-s) \cap b^{-1}([r_1,\infty)))}
 {r^p-s^p} } \ge \\
\qquad &\ge & \frac{Vol(\Omega(K)\cap b^{-1}(r_0,s))}{s^p-r_0^p}
  \,\,   +  \,\, \frac{Vol_{n-1}(\Omega(\bar{U})\cap b^{-1}((s)))}{ps^{p-1}}.
\end{eqnarray}
However, by the strongly minimal volume growth, we have
\begin{eqnarray}
\lim_{r\to\infty}\frac{Vol(B_{x_0}(r+r_2-s) \cap b^{-1}([r_1,\infty)))}
{r^p-s^p} 
&=&\frac {Vol_{n-1}(B_{x_0}(R_0) \cap b^{-1}(r_1))}{p r_1^p} \nonumber\\
&\le& \frac{Vol(\Omega(K)\cap b^{-1}(r_1,s))}{s^p-r_1^p},
\end{eqnarray}
and so, 
$$
Vol_{n-1}(\Omega(\bar{U})\cap b^{-1}((s)))=0,
$$
which contradicts (\ref{VolSlice}).

Thus
$$
b^{-1}([r_1,\infty)) \subset \Omega(K).
$$
This implies that $b^{-1}(r) \subset B_{x_0}(r-r_1+R_0)$ for any
$r \ge r_1$, so it is compact.  

Furthermore, for $r \le r_1$,
$b^{-1}(r)$ is a subset of the closed tubular neighborhood 
$T_{r_1-r}(b^{-1}(r_1))$ as can be seen by using Busemann rays to travel
from $b^{-1}(r)$ to $b^{-1}(r_1)$.  Thus $b$ is a proper function. 
\ProofEnd

\begin{lemma} \label{ConstRatio}
Given $M^n$ as described above, then the ratio
\begin{equation}
V(r)=\frac {Vol_{n-1}(\Omega(K)\cap b^{-1}(r))}{r^{(\frac 1 2 -v)(n-1)}}
\end{equation}
is a constant function for all $r \ge r_1$,
and, thus,
\begin{equation}
\frac{Vol_{n-1}(b^{-1}([r_2,r_3]))}{(r_3^p-r_2^p)}=
\frac{Vol_{n-1}(b^{-1}(r_1)\cap \bar{B}_{x_0}(R_0))}{p r_1^{p-1}}
\end{equation}
for all $r_3 > r_2 \ge r_1$.
\end{lemma}

\noindent
{\bf Proof:}  
Let $K=\bar{B}_{x_0}(R_0) \cap b^{-1}(r_1)$
as in Lemma~\ref{EndinOm}.
By the definition of $\Omega(K)$, we know that for all 
$x \in \Omega(K) \cap b^{-1}(r_1,r)$ there exists 
$x' \in K \subset \bar{B}_{x_0}(R_0)$
such that $x \in \gamma_{x'}([r_1,r])$; so 
\begin{equation} \label{choiK1}
\Omega(K) \cap b^{-1}((r_1,r)) \subset \bar{B}_{x_0}(r-r_1+R_0).
\end{equation}
Note also that
\begin{equation} \label{choiK2}
K \cap b^{-1}(r_1)= \Omega(K)\cap b^{-1}(r_1),
\end{equation}
for our choice of $K$.
Strongly minimal volume growth,
Theorem~\ref{Monotonicity}, (\ref{choiK1}) and (\ref{choiK2})
imply that for any $b>a \ge r_1$
\begin{eqnarray}
\frac{Vol_{n-1}(K \cap b^{-1}(r_1))}{p r_1^{p-1}}
&=& \lim_{r\to\infty}\frac {Vol(B_{x_0}(r-r_1+R_0))\cap b^{-1}([r_1,\infty))}
{r^p-r_1^p}\\
&\ge&\lim_{r\to\infty} \frac {Vol(\Omega(K) \cap b^{-1}(r_1,r))}{r^p-r_1^p}
\label{choiK3}\\
&\ge&\frac {Vol_{n-1}(K \cap b^{-1}(r_1))}{pr_1^{p-1}}.
\end{eqnarray}
Since all the inequalities must be equalities
and the limit in (\ref{choiK3}) is monotone, we get
\begin{equation} \label{VolumeForced1}
\frac{Vol_{n-1}(K\cap b^{-1}(r_1))}{pr_1^{p-1}}=
\frac{Vol(\Omega(K) \cap b^{-1}(r_1,r))}{r^p-r_1^p},
\end{equation}
for any $r > r_1$.
Subtracting (\ref{VolumeForced1}) with $r=a$
from (\ref{VolumeForced1}) with $r=b$, and reworking the equation gives
\begin{equation} \label{VolumeForced2}
\frac {Vol(\Omega(K)\cap b^{-1}(a,b))}{b^p-a^p}
\end{equation}
is a constant with respect to $a$ and $b$.  Then applying
Lemma~\ref{EndinOm}, we obtain the lemma.
\ProofEnd

The fact that the Busemann function is proper in a manifold with
strongly minimal volume growth makes it much easier to prove that
the metric is rigid than if this were not the case.  We need not trace 
through the rather involved
proof of Theorem~\ref{Monotonicity} with its unions of star-shaped sets
about points 
in $b^{-1}(R)$ [So2].  Instead we use the following very simple 
sets.

First we fix $r_3 >r_1 \ge r_0$.  We will prove that
$b^{-1}([r_1, r_3])$ is isometric to the appropriate warped
product manifold $b^{-1}(r_1) \times_{(b/r_1)^{(1/2 -v)}} [r_1,r_3]$
and that $b$ is smooth on $b^{-1}([r_1, r_3])$.  To do so we fix
$\vare_0 < (r_3-r_1)/10$ and take $r_4 = 2 r_3 $.

\begin{definition} \label{InnerAnn}
{\em 
Let $R > r_4$ and $r_1 \le a < b \le R$ we define the {\em inner annulus},
\begin{equation}
S^{r_4}_{a, b, R}\,\,=
\,\, Ann_{R-b,R-a}(\gamma(R)) \,\cap\, S^{r_4}_{\gamma(R)}
\end{equation}
where
$$
S^{r_4}_{\gamma(R)}= \{ \sigma([0,L]): \sigma \textrm{ is a min geod s.t. }
 \sigma(0)=\gamma(R), \sigma(L)\in b^{-1}((-\infty, r_4])\}.
$$
We will often write only $S_{a, b, R}=S^{r_4}_{a, b, R}$.
}
\end{definition}

Recall that $b(x)=\lim_{R\to \infty}(R-\rho_R(x))$
where the convergence is uniform on compact sets.  
So given any
$\vare >0$ and any compact set like $b^{-1}([r_1,r_4])$, there 
exists $R_{\vare}=R^{r_1,r_4}_{\vare}$ such that
\begin{equation} \label {Rvare}
b(x) \le R-\rho_R(x) +\vare \qquad \forall R \ge R_{\vare}.
\end{equation}
On the other hand, by the definition of $b(x)$, we have
\begin{equation} \label{Rvare2}
b(x) \ge \lim_{s \to \infty} \bigg( s -d(x, \gamma(R))-d(\gamma(R), \gamma(s))
\bigg)= R-d(x, \gamma(R)).
\end{equation}
Thus the inner annuli, $S_{a,b,R}$, are close to the compact regions
$b^{-1}([a,b])$ as described precisely in the following lemma.

\begin{lemma} \label{BuseInnerAnn}
Let $r_1 \le a < b \le r_4$.  Then for all $\vare >0$ there
exists $R_{\vare}=R^{r_1,r_4}_{\vare}$ as defined in (\ref{Rvare})
such that
$$
S_{a+\vare, b-\vare, R} \subset b^{-1}([a,b]) \subset S_{a-\vare, b+\vare, R}
$$
for all $R \ge R_\vare$.
\end{lemma}

\noindent{\bf Proof:}
Fix $R \ge R_\vare$.

Suppose $x\in b^{-1}([a,b])\subset b^{-1}([r_1,r_4])$.  Then
by (\ref{Rvare}) and (\ref{Rvare2}), we have
$$
L_x=\rho_R(x) \in [R-b-\vare, R-a+\vare].
$$
Thus there exists a minimal geodesic $\sigma$ from $\sigma(0)=\gamma(R)$
to $\sigma(L_x) =x \in b^{-1}([r_1,r_4])$.  So
$$
x \in S_{a-\vare, b+\vare, R}.
$$

On the other hand, if $y \in S_{a+\vare, b-\vare, R}$, then there exists $L>0$
and a minimal geodesic, $\sigma$, from $\sigma(0)=\gamma(R)$ to 
$\sigma(L) \in b^{-1}([r_1,r_4])$ and there exists 
$t\in [R-b+\vare, \min\{R-a-\vare,L\}]$ such that $y=\sigma(t)$.
Since $\sigma(L)\in b^{-1}([r_1,r_4])$, we can apply (\ref{Rvare})
to get,
$$
b(\sigma(L)) \le R-L+\vare
$$
Thus
$$
b(y) \le b(\sigma(L))+ d(\sigma(L), y) 
 =b(\sigma(L))+ L-d(y, \gamma(R)) \le R +\vare -(R-b+\vare)=b.
$$ 
Meanwhile, by (\ref{Rvare2}) we have,
$$
b(y) \ge R-d(y, \gamma(R)) \ge R-(R-a-\vare) > a
$$
and we have $y \in b^{-1}([a,b])$.
\ProofEnd

Thus the volumes of the regions between Busemann level sets can be compared
to the volumes of these inner annuli.  On the other hand, the volumes of
the inner annuli can be controlled using the Relative Volume Comparison
Theorem on the star-shaped sets, $S_{r_1,R,R}$ about $\gamma(R)$
[BiCr, GrLaPa, Ch].
We can bound the Ricci curvature from below using the following
lemma which is a direct consequence  of the techniques used in the last proof.

\begin{lemma} \label{RicciInnerAnn}
Given any $\vare$, let $\vare_v=\vare$ if $v<1/2$ and
$\vare_v=0$ otherwise.

There exists $R^{r_1,r_4}_\vare$ 
as defined in (\ref{Rvare}) such that for all
$x$ in the star shaped set, $S^{r_1,r_4}_{r_1,R,R}$,
we have
\begin{equation}
Ric_y \ge \frac {(n-1)(\frac 1 4 -v^2)}{(R-d(y,\gamma(R))+\vare_v)^2}
\end{equation}
for all $R \ge R_\vare$.
\end{lemma}

Thus we have the same lower Ricci curvature bound as the comparison 
manifolds of Definition~\ref{CompMan} if we consider
$\gamma(R)$ to be the base point.

We will now control the mean curvature of the level sets of $\rho_R$
within $S_{r_1,R,R}$ using the comparison manifolds, $M^n_{R,\vare}$, of
Definition~\ref{CompMan} combined with the above lemma regarding
Ricci curvature and the Relative Volume Comparison Theorem.  To control
the mean curvature from below we will employ the volume estimates on the
regions between Busemann levels and the relationship between those
regions and the inner annuli.  We begin with some facts and definitions
regarding the mean curvature of a distance function.

\begin{definition} \label{meancurvdefn} {\em
Let $m_R(x)$ be the {\em mean curvature 
with respect to the inward normal } of 
the sphere of radius $d(x, \gamma(R))$ about $\gamma(R)$ evaluated at $x$.
Recall that this mean curvature is well defined along minimal geodesics
to $\gamma(R)$ and is equal to $\Delta \rho(x)$ where it is defined.  See,
for example, [Ch].  

Given $x \in M^n$,
let $\bar{m}_{R,\vare}(x)$ be the mean curvature
with respect to the inward normal of the 
distance sphere of radius $\rho_R(x)$ about the origin in the comparison
warped product manifold, $M^n_{R,\vare}$. } 
\end{definition}

\begin{lemma} \label{meancurv}
The comparison mean curvature satisfies
\begin{equation} \label{meancurv1}
\bar{m}_{R,\vare}(x) =\frac {-(n-1)}{R+\vare-\rho_R(x)}
    \bigg(\frac 1 2  + 
    v\, coth\bigg(Log\bigg(\frac{R+\vare-\rho_R(x)}{R+\vare}\bigg) \bigg) 
    \bigg)
\end{equation}
and as $R$ goes to infinity, we have
\begin{equation} \label{meancurv2}
\lim_{R\to\infty}\bar{m}_{R,\vare}(x)=
    \frac {-(n-1)}{b(x)}\bigg( \frac 1 2 + v (-1)\bigg)
\end{equation}
uniformly on compact sets.
\end{lemma}

The proof of the lemma is an exercise on warped product manifolds and
can be found in [So1].

We now need to show that the mean curvature of spheres around $\gamma(R)$,
$m_R(x)$ approaches $\bar{m}_{\vare, R}(x)$ as $R$ approaches infinity and
$\vare$ goes to $0$.  Then, in some weak sense we will have
\begin{equation} \label{roughidea}
\Delta b(x)=-\lim_{R\to \infty} \Delta \rho(x)
= -\lim_{R\to \infty} m_{R}(x)   
= \frac {-(n-1)(1/2-v)}{b(x)}
\end{equation}
Since $\Delta \rho_R(x)$ and $\Delta b(x)$ are defined on different
subsets of M, we need to prove (\ref{roughidea}) very carefully.  
In the next lemma we will only obtain a weak estimate on 
$\lim_{R\to \infty} m_R(x)$ using the strongly minimal volume
growth, but that will suffice.

\begin {lemma} \label{DiffA_R(t)}
Let $M^n$ have the properties defined at the beginning of this section.
Fix $r_3>r_1$.
Then for all $\delta >0$, there exists $\vare_\delta >0$ 
such that for all $\vare <\vare_\delta$ there exists
$R_\delta(\vare,r_1,r_3)>2r_3$ 
and there exists a constant, $C_{r_1,r_3,R+\vare}$, such that
$\lim_{R\to\infty} C_{r_1,r_3,R+\vare}=C_{r_1,r_3}>0$
such that for all $R \ge R_\delta(\vare,r_1,r_3)$, we have,
\begin{equation}\label{DiffA_R(t)2}
0<\int_{s=r_1}^{r_3}\int_{S_{s,r_3,R}(K)}
   \bigg(\Delta \rho_R(x)-\bar{m}_{R,\vare}(x)\bigg)
           \,dvol\,dt < \frac{\delta}{ C_{r_1,r_3,R+\vare}}.
\end{equation}
Furthermore, the integrand is nonnegative.
\end{lemma}

In this lemma $dvol$ is the volume form on $M^n$.  The integral over
$S_{a,b,R}(K)$ can be best understood using the following definition.
Recall Definition~\ref{InnerAnn}.

\begin{definition} \label{Theta}
Let $\Theta(R) \subset TM_{\gamma(R)}$ be defined as
$$
\Theta_R=\Theta^{r_4}_R= \{\sigma'(0):\sigma \textrm{ is a min geod s.t. }
 \sigma(0)=\gamma(R), \sigma(L)\in b^{-1}((-\infty, r_4])\},
$$

Given any $\theta\in\Theta_R$, let $A_R(t,\theta)$ be
the warped product of the Jacobi fields along
$exp_{\gamma(R)}t\theta$ which are 0 when $t=0$ and whose first covariant
derivatives are orthonormal when $t=0$.  We set $A_R(t,\theta)$ to be 
continuous up to and including the cut point.  After a cut point it is set to
be 0.  

So $dvol=A_R(t,\theta) \,\nu_\theta\, dt$ where
$\nu_\theta$ is the volume form on the unit sphere.  

\end{definition}

\noindent
{\bf Proof of Lemma~\ref{DiffA_R(t)}:}  
Let $r_4 = 2 r_3 > r_3$.  Let $S_{a,b,R}=S^{r_4}_{a,b,R}$ as in 
Definition~\ref{InnerAnn}.

Given any $\vare >0$ and $R>r_3$, let 
\begin{equation}
J(t)=J_{\vare_v,R}(t)
\end{equation}
be the Jacobi field defined in Definition~\ref{CompMan}
where $\vare_v=\vare$ if $v<1/2$ and
$\vare_v=0$ otherwise.
By the Bishop Volume Comparison Theorem [BiCr],
the lower bound on Ricci curvature in the comparison manifolds,
(\ref{RicciComp1}), and the Ricci curvature bound on $S_{a,b,R}$,
[Lemma~\ref{RicciInnerAnn}],  we know that 
$\frac { A_R(t,\theta) }{ J(t)^{n-1} }$ is nonincreasing
for $t \in [0, R-r_1]$ and $\theta \in \Theta_R$.  
We wish to show that this ratio is almost constant.

First we note that, by Definition~\ref{Theta} and (\ref{Rvare2}),
if $a \le b \le r_4$ and $\theta \in \Theta_R$
such that $exp_{\gamma(R)}(t \theta)$ is minimal for $t\in [0,R-a]$,
then 
$$
b(exp_{\gamma(R)}((R-a) \theta)) \le R-(R-a) \le r_4.
$$
So by Definition~\ref{InnerAnn},
$exp_{\gamma(R)}(t \theta)\in S_{a,b,R}$ for all $t\in [R-b, R-a]$
and we have
\begin{equation}
Vol(S_{a,b,R}) = \int_{R-b}^{R-a}\int_{\Theta_R}A_R(t,\theta)
\,d\nu_{\theta}\,dt
\end{equation}

Thus, by the Relative Volume Comparison Theorem,
\begin{equation} \label{StoAtoJ1}
\frac{Vol(S_{r_3,r_4, R})}{Vol(S_{r_1,r_3, R})} = 
\frac{\int_{R-r_4}^{R-r_3} \int_{\theta_{R}} A_R(t,\theta) \,d\nu_\theta \,dt }
{\int_{R-r_3}^{R-r_1} \int_{\theta_{R}} A_R(t,\theta) \,d\nu_\theta \,dt }
\ge 
\frac{\int_{R-r_4}^{R-r_3}J_{R.\vare}(t)^{n-1}  \,dt }
{\int_{R-r_3}^{R-r_1} J_{R.\vare}(t)^{n-1} \,dt }.
\end{equation}
By Lemma~\ref{BuseInnerAnn}, 
we know that for any $\vare >0$, there exists $R_\vare$ such that
\begin{equation} \label{StoAtoJ2}
\frac{Vol(S_{r_3,r_4, R})}{Vol(S_{r_1,r_3, R})}  
\le \frac {Vol(b^{-1}([r_3-\vare, r_4+\vare]))}
{Vol(b^{-1}([r_1+\vare, r_3-\vare]))}\qquad \forall R \ge R_\vare.
\end{equation}

By the strongly minimal volume growth and Lemma~\ref{ConstRatio}, we have
\begin{equation} \label{StoAtoJ3}
\frac {Vol(b^{-1}([r_3-\vare, r_4+\vare]))}
{Vol(b^{-1}([r_1+\vare, r_3-\vare]))}
=\frac {(r_4+\vare)^p-(r_3-\vare)^p}{(r_3-\vare)^p-(r_1+\vare)^p}
\end{equation}
On the other hand, by (\ref{JacobiLimit}),
\begin{equation} \label{StoAtoJ4}
\lim_{R\to\infty}\frac{\int_{R-r_4}^{R-r_3}J_{R.\vare}(t)^{n-1}  \,dt }
{\int_{R-r_3}^{R-r_1} J_{R.\vare}(t)^{n-1} \,dt }
=\frac {(r_4)^p-(r_3)^p}{(r_3)^p-(r_1)^p}
\end{equation}

Given any $\delta' >0$, there exists $\vare_{\delta'}>0$ small enough that
the right hand side of (\ref{StoAtoJ4}) is within $\delta'$ of
the right hand side of (\ref{StoAtoJ3}) for all $\vare < \vare_{\delta'}$, 
and there exists $\bar{R}_{\delta'}$ large enough that the 
the right hand side of (\ref{StoAtoJ4}) is within $\delta'$ of
the right hand side of (\ref{StoAtoJ1}).  Thus the inequalities
in (\ref{StoAtoJ1})-(\ref{StoAtoJ4}) are almost equalities,
and we have the following.

For all $\delta >0$, let $\delta'=\delta/(2Vol(S_{r_1,r_3,R}))$,
there exists $\vare_{\delta'}>0$, such that
for all $\vare < \vare_{\delta'}$,
there exists
\begin{equation} \label{StoAdefR}
R_{\delta}(\vare,r_1,r_3) = \max \{ R_{\vare}, \bar{R}_{\delta'} \},
\end{equation}
such that we have
\begin{equation} \label{StoAtoJ5}
0 \le 
\frac{\int_{R-r_4}^{R-r_3} \int_{\theta_{R}} A_R(t,\theta) \,d\nu_\theta \,dt }
{\int_{R-r_3}^{R-r_1} \int_{\theta_{R}} A_R(t,\theta) \,d\nu_\theta \,dt }
-
\frac{\int_{R-r_4}^{R-r_3}J_{R.\vare}(t)^{n-1}  \,dt }
{\int_{R-r_3}^{R-r_1} J_{R.\vare}(t)^{n-1} \,dt }.
\le 2 \delta'.
\end{equation}

By multiplying both sides of (\ref{StoAtoJ5}) by 
$Vol(S_{r_1,r_3,R})$ and multiplying both sides by
$\int_{R-r_3}^{R-r_1} J_{R.\vare}(t)^{n-1} \,dt$, we have,
\begin{eqnarray}
0 &\le& 
\int_{R-r_4}^{R-r_3} \int_{\theta_{R}} A_R(t,\theta)\,d\nu_\theta \,dt
\int_{R-r_3}^{R-r_1}J(t)^{n-1}dt\,\, \nonumber   \\
& &   -\,\,
\int_{R-r_3}^{R-r_1} \int_{\theta_{R}} A_R(t,\theta) \,d\nu_\theta \,dt 
\int_{R-r_4}^{R-r_3}J(t)^{n-1}dt   \nonumber       \\
\label{CompRay}
& \le &
\delta \,\,\,\int_{R-r_4}^{R-r_3}\,\,\,(\,J(t)\,)^{n-1}\,\,dt.
\end{eqnarray}

As mentioned above,
$\frac {A_R(t,\theta)}{J(t)^{n-1}}$ is nonincreasing.
We wish to show that 
this ratio is almost constant using (\ref{CompRay}).

Note that,
\begin{equation}
\int_{R-b}^{R-a}\int_{\Theta_R}A_R(t,\theta)
\,d\nu_{\theta}\,dt\\
=\int_{R-r}^{R-r_3}\int_{\Theta_R}
\frac{A_R(t,\theta)}{J(t)^{n-1}}J(t)^{n-1}
\,d\nu_{\theta}dt.
\end{equation}

By substituting this into (\ref{CompRay})
and by subtracting and adding
$$
\int_{\Theta_R}\frac{A_R(R-r_3,\theta)}
{J(R-r_3)^{n-1}}\,d\nu_{\theta}\,\,
\int_{R-r_3}^{R-r_1}J(t)^{n-1}dt\,\,
\int_{R-r}^{R-r_3}J(t)^{n-1}dt
$$ 
to (\ref{CompRay}), we get
$$
\delta \int_{R-r_4}^{R-r_3}J(t)^{n-1}dt 
>I_1 \int_{R-r_3}^{R-r_1}J(t)^{n-1}dt
+I_2\int_{R-r_4}^{R-r_3}J(t)^{n-1}dt,
$$
where
\begin{equation}
I_1=\int_{R-r_4}^{R-r_3}\int_{\Theta_R}\left(\frac{A_R(t,\theta)}
{J(t)^{n-1}}
-\frac{A_R(R-r_3,\theta)}{J(R-r_3)^{n-1}}\right)J(t)^{n-1}
\,d\nu_{\theta}\,dt \label{Ione}
\end{equation}
and
\begin{equation}
I_2=\,\,\int_{R-r_3}^{R-r_1}\int_{\Theta_R}
\left(\frac{-A_R(t,\theta)}{J(t)^{n-1}}+
    \frac{A_R(R-r_3,\theta))}
{J(R-r_3)^{n-1}}\right)J(t)^{n-1}
d\nu_{\theta}dt \label{Itwo}.
\end{equation}

Since $\frac {A_R(t)}{J(t)^{n-1}}$ is decreasing, the integrands of
$I_1$ and $I_2$ are 
both positive.  Thus, $I_2 < \delta$ and
\begin{equation} \label{BG2}   
0<\int_{R-r_3}^{R-r_1} \int_{\Theta_R} \bigg|
\frac{A_R(R-r_3,\theta)}{(J(R-r_3))^{n-1}}-
                \frac{A_R(t,\theta)}{(J(t))^{n-1}}\bigg|
J(t)^{n-1}\,d\nu_{\theta}\,dt < \delta.
\end{equation}
Thus the ratio, $\frac {A_R(t,\theta)}{J(t)^{n-1}}$, is almost constant. 

Given $\theta \in \Theta_R$, let $d_\theta \in (R-r_4, \infty]$ 
be the distance to the
cutpoint of $\gamma(R)$ along $exp_{\gamma(R)}(t \theta)$.
Then $A_R(t,\theta)$ is smooth for $t \in (0,d_{\theta})$,
continuous for $t \in [0,d_\theta]$, and $0$ for $t>d_\theta$.  
Since the mean curvatures are evaluated with respect to the inward
normal, we also know that
\begin{equation} \label{Delta}
\frac{A'_R(t,\theta)}{A_R(t,\theta)}
=-m_R(exp_\gamma(R)(t \theta))=-\Delta \rho_R\big(exp_\gamma(R)(t \theta)\big),
\end{equation}
for $t\in (0,d_\theta)$, where 
$A'_R(t,\theta)=\frac {\partial}{\partial t}A_R(t,\theta)$.
See, for example [Ch].

Let $r_{\theta,t}=max\{min\{t,d_\theta\},R-r_3\}$.  
We will use this function to seperate off the 
differentiable part of the integrand of (\ref{BG2}).
$$
0< \int_{R-r_3}^{R-r_1}\int_{\Theta_R} \bigg(
  \frac{A_R(R-r_3,\theta)}{J^{n-1}(R-r_3)}
   -\frac{A_R(r_{\theta,t},\theta)}{J^{n-1}(r_{\theta,t})} 
   \bigg)J^{n-1}(r_{\theta,t})\,d\nu_\theta\,dt \,\,\,\,\,\,
$$
$$
\,\,\,\,\,\,\,\,   +\,\,\,\, \int_{R-r_3}^{R-r_1}\int_{\Theta_R} \bigg(
\frac{A_R(r_{\theta,t},\theta)}{J^{n-1}(t)}
   -\frac{A_R(t,\theta)}{J^{n-1}(t)}
    \bigg)J^{n-1}(r_{\theta,t})\,d\nu_\theta\,dt < \delta.
$$

Since $t\ge r_{\theta,t}\ge R-r_3$ and the ratios are decreasing, 
both integrands are positive.  Thus we have,
$$
0< \int_{t=R-r_3}^{R-r_1}\int_{\Theta_R} \bigg(
  \frac{A_R(R-r_3,\theta)}{J^{n-1}(R-r_3)}
   -\frac{A_R(r_{\theta,t},\theta)}{J^{n-1}(t)} 
   \bigg)J^{n-1}(t)\,d\nu_\theta\,dt < \delta.
$$
The function $A_R$ is differentiable in this integral, so
$$
0< \int_{t=R-r_3}^{R-r_1}\int_{\Theta_R} \bigg(
  -\int_{l=R-r_3}^{r_{\theta,t}}
   \frac{d}{dl}\bigg(\frac{A_R(l,\theta)}{J^{n-1}(l)}\bigg) 
    dl \bigg) J^{n-1}(t)\, d\nu_\theta\, dt  < \delta,
$$
and
$$ 
-\int_{t=R-r_3}^{R-r_1}\int_{\Theta_R} \int_{l=R-r_3}^{r_{\theta,t}} 
  \bigg(\frac{A'_R(l,\theta)}{A_R(l,\theta)}-\frac{(J^{n-1}(l))'}{J^{n-1}(l)}
  \bigg)\frac {J^{n-1}(t)}{J^{n-1}(l)} A_R(l,\theta)
   dl   d\nu_\theta  dt  < \delta.
$$

Recall that $S_{R-t,r_3,R}$ only consists of points on minimizing 
geodesics from $\gamma(R)$ to $b^{-1}((-\infty, r_4])$.  Thus
$$
S_{R-t,r_3,R}=\{exp_{\gamma(R)}(l\theta): \theta \in \Theta_R, 
l\in [R-r_3, r_{\theta,t}) \}
$$ 
and $exp_{\gamma(R)}$ is invertible on $S_{R-t,r_3,R}$.
Let $l=\rho_R(x)$ and $\theta$ be defined such that 
$x=exp_{\gamma(R)}(l\theta)$ where $x\in S_{R-t,r_3,R}$.
Furthermore $dvol=A_R(l,\theta)\,dl \, d\nu_\theta$ on $S_{R-t,r_3,R}$.

So we have
$$
0<-\int_{t=R-r_3}^{R-r_1}\int_{S_{R-t,r_3,R}}
  \bigg(\frac{A'_R(l,\theta)}{A_R(l,\theta)}-\frac{(J^{n-1}(l)'}{J^{n-1}(l)}
     \bigg) \,dvol\,dt < \frac{\delta}{ C_{r_1,r_3,R+\vare_v}}
$$
where 
\begin{equation} \label {EqnChere}
C_{r_1,r_3,R+\vare_v}=\min_{(R-r_1\ge t \ge l\ge R-r_3)}\bigg(
            \frac {J_{R,{\vare_v}}^{n-1}(t)}{J_{R,{\vare_v}}^{n-1}(l)}\bigg).
\end{equation}

Since this integral avoids cut points of $\gamma(R)$, we have
$$
0<-\int_{t=R-r_3}^{R-r_1}\int_{S_{R-t,r_3,R}}
   \bigg(-\Delta \rho_R(x)+\bar{m}_{R,\vare}(x)\bigg)
           \,dvol\,dt < \frac{\delta}{ C_{r_1,r_3,R+\vare}}.
$$
Setting $s=R-t$, we are done proving (\ref{DiffA_R(t)2})
with $R_\delta(\vare,r_1,r_3)$ defined in (\ref{StoAdefR}).

Finally, for all $v\in [0,\frac{n+1}{2(n-1)}]$, we know that
\begin{eqnarray*}
\lim_{R\to\infty}C_{r_1,r_3,R+\vare} &= &\lim_{R\to\infty}
    \min_{r_1\le t'\le l'\le r_3}\frac
   {-(t'+\vare)^{v+1/2}+(R+\vare)^{2v}(t'+\vare)^{-v+1/2}}
      {-(l'+\vare)^{v+1/2}+(R+\vare)^{2v}(l'+\vare)^{-v+1/2}}\\
&=&\min_{r_1\le t'\le l'\le r_3}\frac{(t'+\vare)^{-v+1/2}}{(l'+\vare)^{-v+1/2}}
\,\,\, > \,\,\, 0.
\end{eqnarray*}
\ProofEnd

In Lemma~\ref{DiffA_R(t)}, the estimate on $\Delta \rho_R(x)$ involves
a double integral.  This can be studied as a single integral over
a region in the isometric product, $M \times \RR$, where we define 
the {\em extended Busemann function}, $\bar{b}(x,s)=b(x)$, and the 
{\em extended distance function}, $\bar{\rho}_R(x,s)=\rho(x)$.
Note that the Laplacians of these extended
functions on $M \times \RR$ are the same as the original functions' 
Laplacians in $M$ because the extended functions are constant in $s$.

\begin{proposition}\label{=CompRay2}
Let $M^n$ be as defined above.

Then the Busemann function with respect to that ray, $b(x)$, is smooth and
\begin{equation} \label{DiffEqCR2}
\Delta b(x) = \frac {(n-1)(1/2-v)}{b(x)}
\end{equation}
on $b^{-1}((r_1,\infty))$. 
\end{proposition}

\noindent
{\bf Proof:}

Choose any $x_1 \in b^{-1}((r_1, \infty))$.  Let 
$h =\frac {b(x_1)-r_1}{4}$, and let
$$
b_1=b(x_1)-h , 
b_2=b(x_1)+h,  \textrm{ and } r_3=b(x_1)+ 2h.
$$

Define the open set, $U$, as follows: 
\begin{equation} \label{DefnU}
U = \{(x,s):\,\, s \in (b_1, b_2), \,\,x \in b^{-1}(s, b_2) \}
\subset M \times \RR.
\end{equation}
We will prove that 
$$
\Delta \bar{b}(x,s) =\frac {(n-1)(1/2-v)}{\bar{b}(x,s)}
$$
in the weak sense on $U$.

By Lemma~\ref{BuseInnerAnn}, there exists $R_h^{r_1,r_3}$ such that
$$
b^{-1}([s,b_2]) \subset S_{s-h, b_2 +h, R}^{2r_3}
\subset S_{s-h, r_3, R}^{2r_3}.
$$
for all $R \ge R_h^{r_1,r_3}$.  Thus, since $b_1-h>r_1$,
\begin{equation} \label{UinExtF}
U \subset \{(x,s):\, s \in (r_1,r_3),\, x \in S_{s, r_3, R}^{2r_3}\}.
\end{equation}

Let $\phi:U \to \RR$ be a smooth nonnegative function with compact support.
By Lemma~\ref{DiffA_R(t)} and (\ref{UinExtF}), we know that for all $\delta>0$
there exists $\vare_\delta>0$, such that for all $\vare <\vare_\delta$,
there exists $R_\delta(\vare, r_1,r_3)$ and $C_{r_1,r_2,R+\vare}$
as defined in Lemma~\ref{DiffA_R(t)}, such that 
\begin{equation} \label{cutoff1}
0<\int_U \,\phi(x,s)
   \bigg(\Delta \rho_R(x)- \bar{m}_{R,\vare}(x)\bigg)
           \,dvol\,dt < \frac{\delta}{ C_{r_0,r_1,R+\vare}}
\end{equation}
Here we have used the fact 
that the integrand of (\ref{DiffA_R(t)2}) is nonnegative
and the integrand here is still nonnegative.

Note that $\Delta \rho_R(x)=-\Delta (R-\rho_R(x))$ 
since $R$ is just a constant.  Note also that $\rho_R(x)=\bar{\rho}_R(x,s)$.
Using the fact that the cut off function, $\phi$, is $0$ near the boundary
to integrate by parts, we get
\begin{eqnarray}
\lefteqn{\int_U \,\phi(x,s)
 \bigg(\Delta \rho_R(x) \bigg)\,dvol\,dt=} \nonumber\\
&=&\int_U \,\phi(x,s)
 \bigg(-\Delta(R-\bar{\rho}_R(x,s)) \bigg)\,dvol\,dt\nonumber\\
&=&\int_U\,
- \Delta\phi(x,s)(R-\bar{\rho}_R(x,s))\,dvol\,dt \label{GenLapHere}
\end{eqnarray}
By the uniform convergence of $b(x)$ on $U \subset b^{-1}([r_1,r_3])$,
we know there exists $R_{\delta,r_1,r_3}$ such that
\begin{equation}\label{ConvBarB}
|R-\bar{\rho}_R(x,s) -\bar{b}(x,s)| < \delta \quad
\forall R \ge R_{\delta,r_1,r_3}.
\end{equation}

Furthermore, by Lemma~\ref{meancurv}, there exists
$R_\delta$ such that 
\begin{equation}\label{meancurvhere}
\Bigg|\bar{m}_{R,\vare}(x)-
    \frac {-(n-1)}{b(x)}\bigg( \frac 1 2 + v (-1)\bigg)\Bigg|
< \delta \quad
\forall R \ge R_{\delta}
\end{equation}

Substituting (\ref{GenLapHere}),(\ref{ConvBarB}) and (\ref{meancurvhere})
all into (\ref{cutoff1}), for 
$$
R > \max \{R_\delta,R_{\delta_0,r_1,r_3},R_\delta(r_1,r_3)\},
$$ 
we have 
$$
\int_U
   -\Delta\phi(x,s)\bar{b}(x,s)
       + \phi(x,s)\bigg(\frac {(n-1)(\frac 1 2 -v)}{\bar{b}(x,s)}\bigg)
           \,dvol\,dt < \frac{\delta}{C_{r_0,r_1}} + 2 \delta Vol(U)
$$
where the integrand is nonnegative.
This equation no longer depends on $R$ or $\vare$, so it holds
for all $\delta >0$.  Taking $\delta$ to $0$ we see that 
$$
\Delta \bar{b}(x,s)=\frac {(n-1)(1/2-v)}{\bar{b}(x,s)}
$$
in the generalized sense on $U$. 

By elliptic regularity, we know that $\bar{b}(x,s)$ is differentiable 
on $U$ and satisfies
\begin{equation} \label{laplbagain}
\Delta \bar{b}(x,s)= \frac {(n-1)(1/2-v)}{b(x)}.
\end{equation}
The fact that $\bar{b}(x,s)$ appears on both
sides of (\ref{laplbagain}) allows us to pull up its 
differentiability by its bootstraps as high as we want.

Note that
$$
B_{x_1}(h/2) \times (b_1, b_1 +h/2) \subset U,
$$
so $b(x)=\bar{b}(x,b_1 +h/4)$ is smooth on $B_{x_1}(\vare/2)$
and satisfies the differential equation, (\ref{DiffEqCR2}).
This can be done at each $x_1 \in b^{-1}((r_1,\infty))$, so
we have proven the proposition.
\ProofEnd

We can now use this proposition combined with the Bochner Weitzenboch 
formula to prove that $b^{-1}((r_1,\infty))$ is a warped product
and thus the Busemann function is smooth on $b^{-1}([r_1,\infty))$.

\begin{lemma}\label{BochnerWarp}
Let $M^n$ satisfy (\ref{RicciBound}) with $v\in [0, (n+1)/2(n-1))$.
Suppose the Busemann function, $b(x)$, is a solution
of
\begin{equation}\label{LapB5}
\Delta b(x) = \frac {(n-1)}{b(x)}\bigg( \frac 1 2 -(v)\bigg)
\end{equation}
on some subset, $b^{-1}((r_1,\infty)) \subset b^{-1}([r_0,\infty))$,

Then $b^{-1}([r_1,\infty))$ is isomorphic to the warped product,
$$
b^{-1}([r_1,\infty))=
b^{-1}(r_1) \, \times_{(b/r_1)^{(1/2 -v)}} \, [r_1, \infty).
$$
\end{lemma}

\noindent
{\bf Proof:}
First, we write the Bochner Weitzenboch Formula applied to $\grad b$.
That is, we
substitute $\grad b = \xi$ into the formula,
$$
\frac 1 2 \Delta (\xi^i\xi_i)=g^{rs}g^{ab}\xi_{a,r}\xi_{b,s} + 
g^{rs}g^{ab}(\xi_{a,rs}+\xi_{r,as}-\xi_{r,sa})\xi_b + R_{lm}\xi_l\xi_m,
$$
of [Boc, Lemma 2], to get
$$
\frac 1 {2} \Delta |\grad b |^2=
        |Hess\,\, b|^2 + <\grad \Delta b,\grad b> + <Ric \grad b, \grad b>.
$$
We  now use the fact that $ |\grad b | =1$, the differential equation
(\ref{LapB5})
and the Ricci bound to get
$$
0 \ge |Hess\,\, b|^2 + (1/2 -(v))(n-1)<\grad b^{-1},\grad b>
                 + (n-1)(1/4 -v^2)/b^2.
$$
Thus,
$$
0 \ge |Hess\,\, b|^2 +(1/2-(v))(n-1)<-b^{-2}\grad b,\grad b>+(n-1)(1/4 -v^2)b^{-2},
$$
and
$$
|Hess\,\, b|^2 \le (n-1)b^{-2}(1/2-v)(1-(1/2 +v))=(n-1)b^{-2}(1/2-v)^2.
$$
On the other hand,
\begin{equation} \label{Cauchy}
|Hess\,\, b|^2=\sum_{i\neq j} b_{i,j}^2 +\sum_{i\ne 1} b_{i,i}^2
    \ge 0 + \frac 1 {n-1} \bigg(\sum_{i\ne 1} b_{i,i}\bigg)^2 
\end{equation}
by the Cauchy-Schwartz inequality and the fact that $b_{1,1}=0$.
So,
$$
|Hess\,\, b|^2\ge \frac 1 {n-1} \big((1/2 -(v))(n-1)b^{-1}\big)^2
   =(n-1)b^{-2}(1/2-v)^2.
$$
Thus the inequalities must be equalities in the Cauchy-Schwartz inequality,
(\ref{Cauchy}), so
$$
b_{i,j}=0 \,\,\,\forall i,j \,\,\,and
\,\,\,b_{k,k}=b_{l,l}\,\,\, \forall k,l \ne 1.
$$
Using the formula for the Laplacian of $b$ once again we get
$$
b_{i,i}=(1/2-v)/b,
$$
and so we can solve for the warping function, $f(b)$,
$$
f''(b)= \frac{(1/2-v) f(b)} { b} \,\,\, \Longrightarrow \,\,\,
      f(b)=\Bigg(\frac b {r_0}\Bigg)^{(\frac 1 2 -v)}f(r_0).
$$
Thus $b^{-1}((r_1,\infty))$ is the desired warped product.  

To complete the proof, note that the boundary,
$b^{-1}(r_1)$, must be isometric to a rescaled $b^{-1}(r)$ for any
$r>r_1$.  So it is smooth and can be included in the warped product. 
\ProofEnd

This completes the proof of Theorem~\ref{VolMinEndWarp}.

Note that if we only had the condition of minimal volume growth,
then we could prove a series of lemmas similar to the ones proven here
to show that $\Delta \bar{b}$ is approximately equal to
$(n-1)(1/2 -v)/ \bar{b}$ in a weak sense.  We could then apply
Cheeger and Colding's Almost Rigidity Theory to prove the
Theorem~\ref{GHWarp}.  Rather than imitating their methods from scratch,
we will adapt one of their key theorems to our situation.  However,
the reader should understand that it is the control on the weak Laplacian of
the Busemann function that gives us the almost rigidity.

\section{Minimal Volume Growth and Almost Rigidity}

\Pind
In this section, we examine the asymptotic properties of a manifold with 
a quadratically decaying lower Ricci curvature bound, (\ref{RicciBound}), and
minimal volume growth.  [Recall Definition~\ref{MinVolGrowth}].
We wish to show that compact regions in such a
manifold are Gromov-Hausdorff close to warped product manifolds 
[Theorem~\ref{GHWarp}].  To do so, we will apply the following theorem
proven in [So2].

\begin{theorem} \label{LinDiam} \label{CompactLevels} [So2, Thm 19]
Let $M$ be a manifold with nonnegative Ricci curvature,
a quadratically decaying lower Ricci curvature bound, (\ref{RicciBound})
with $v\in [0,1/2]$, and minimal volume growth.
Then the Busemann function, $b(x)$, has compact level sets
and their diameter grows at most linearly,
\begin{equation}\label{DiamLin}
diam(b^{-1}(r)) \le C_D |r+1| \qquad \forall r \ge r_0. 
\end{equation}
\end{theorem}

In order to apply this theorem, we will assume that all our manifolds have 
globally nonnegative Ricci curvature for the remainder of the paper.  
Thus the regions $b^{-1}([r,r+L])$ are compact.  These are the compact
regions which are proven to be
close to warped product manifolds in the Gromov-Hausdorff sense
in Theorem~\ref{GHWarp}.

The precise statement of the almost rigidity theorems,
Theorems~\ref{GHWarp} and~\ref{GHCross}  will appear in 2.2 after the 
Gromov-Hausdorff distance and related concepts are defined.

Before going on, it is important to note the following facts from
[So2, Cor 23] reviewed in Section 1.1.  The function,
\begin{equation} 
V(r)=\frac {Vol_{n-1}(b^{-1}(r))}{r^{(1/2 -v)(n-1)}}
\end{equation}
is nondecreasing as a function of $r$ in any manifold with
a quadratically decaying lower Ricci curvature bound, and, if
the manifold also has nonnegative Ricci curvature and minimal
volume growth, then
\begin{equation} \label{Vinfty}
\lim_{r\to\infty}\frac {Vol_{n-1}(b^{-1}(r))}{r^{(1/2 -v)(n-1)}}
=V_\infty<\infty.
\end{equation}  
In fact $V_\infty =V_0$ of the minimal volume growth definition
[Defn~\ref{MinVolGrowth}] unless the
manifold splits isometrically, in which case $V_\infty=V_0/2$.

Both the constants $C_D$ and $V_\infty$ will be refered to in 
the proofs of our almost rigidity theorem and our diameter growth
estimate.

\subsection{Almost Rigidity and the Gromov-Hausdorff Metric}

There are a number of equivalent definitions of the Gromov-Hausdorff
metric on the space of metric spaces.  Here we will use the Gromov-Hausdorff
map to define this metric, since ultimately we will use both
the Gromov-Hausdorff closeness and the particular Gromov-Hausdorff map
to prove our diameter theorem.  See [GrLaPa] for more details.

\begin{definition} \label{DefnGH}
{\em Given $\vare >0$, 
the {\em Gromov Hausdorff distance }, $d_{GH}(X,Y)$ between two compact
metric spaces, $X$ and $Y$, is less than $\vare$ if there
exists a {\em Gromov-Hausdorff map}, $F_{GH}:X \mapsto Y$ 
which is {\em $\vare$-almost onto},
\begin{equation}
T_{\vare}\big(F_{GH}(X)\big)\, \supset \,Y,
\end{equation}
and {\em $\vare$-almost distance preserving},  }
\begin{equation}
\big|\,d_Y(F_{GH}(x_1), F_{GH}(x_2))\,-\,d_X(x_1, x_2) \,\big|\, < \vare.
\end{equation}
The Gromov-Hausdorff map need not be continuous.
\end{definition}

Note that this definition is not quite symmetric.  However, if there exists
$F_{GH}:X \mapsto Y$ with the above properties then the map
$\bar{F}_{GH}:Y \mapsto X$, such that
$\bar{F}_{GH}(y)$ equal any $x\in X$ such that $d_X(y, F_{GH}(x)) <\vare$,
is $2\vare$-almost distance preserving and $2\vare$-almost onto.

\begin{definition} \label{DefnDist}
{\em Given any subset, $U$, of a length space, $N$, we can define a 
{\em localized distance function}, 
\begin{equation}
d_U(x,y)= \inf \{ L(c([0,1])): c(0)=x, c(1)=y, c([0,1])\subset U\}.
\end{equation}

In particular, if we chose any constants $\alpha$ and $ \alpha'$ such that 
$0<\alpha'<\alpha<(b-a)/2$ then we can define 
\begin{equation}
d^{\alpha'}(x,y)=d_{r^{-1}(a+\alpha', b-\alpha')}(x,y)
\end{equation} 
as a distance function on $r^{-1}(a+\alpha', b-\alpha')$ and its restriction
$d^{\alpha,\alpha'}$ to the subset $r^{-1}(a+\alpha, b-\alpha)$.
There is a discussion of these two functions
in [ChCo, Section 3].  }
\end{definition}

\begin{note} \label{components}  
{\em A {\em localized component} of $r^{-1}((a+\alpha, b-\alpha))$ is
a set of the form $U \cap r^{-1}((a+\alpha, b-\alpha))$ where $U$ is
a connected component of $r^{-1}((a+\alpha, b-\alpha))$.  Thus 
$d^{\alpha, \alpha'}(x,y)$ is finite iff $x$ and $y$ are in the same
localized component.

When we say that two spaces are Gromov-Hausdorff close each of which has
more than one such component, then we have paired off all the 
localized components and
shown that each pair is Gromov-Hausdorff close.  In particular, there are the
same number of localized components [ChCo]. 
Note that in a warped product manifold localized components are connected 
components.}
\end{note}

We now present a particular
theorem of Cheeger and Colding which is especially useful 
in the study of manifolds with minimal volume growth. 
This theorem
states that manifolds with lower Ricci curvature bounds and 
{\em almost maximal volume} with respect to a distance function, $r$,
are Gromov-Hausdroff close to certain warped product manifolds 
[ChCo, Thm 4.85]. First we provide a definition of {\em
almost maximal volume}.

Let $m_x(r^{-1}(a))$ denote the mean curvature of $r^{-1}(a)$
at the point $x$.  We will omit the subscript $x$ when it is unimportant.

\begin{definition} \label{AlmostMax}
Let $N^n$ be a Riemannian manifold and $K \subset N^n$ compact.
Let $r(x)=d(x,K)$ be the distance function to $K$.  Fix $b>a\ge 0$.

If the region $r^{-1}(a,b) \subset N^n$ has the following three
properties for some positive smooth function, $f$ and some $\omega\ge0$ 
\begin{equation}
Ric_{N^n}(x) \ge -(n-1) \frac {f''(r(x))}{f(r(x))}  
\end{equation}
\begin{equation}
m(r^{-1}(a))
\le (n-1) \frac{f'(a)}{f(a)}\,\,\,\,on\,\,\,\,r^{-1}(a)  
\end{equation}
\begin{equation}
\frac{Vol(r^{-1}(a,b))}{Vol_{n-1}(r^{-1}(a))}\ge 
(1-\omega)\frac{\int_a^b \, f^{n-1}(r)\,\,dr}{f^{n-1}(a)}   
\end{equation}
then we say that the region $r^{-1}(a,b)$ has $\omega$-almost maximal volume
with respect to the function $f$.
\end{definition}

\begin{theorem} {\bf (Cheeger and Colding)} \label{ChCo4.85}
{\em [ChCo, Thm 4.85]} \newline
Let $N^n$ have $Ricci \ge \Lambda$.  Let $K\subset N^n$ be compact 
and let $r(x)=d(x,K)$.Let $f$ be a smooth
nonnegative function.

Suppose a region $r^{-1}((a,b))\subset N^n$ has 
$\omega$-almost maximal volume with respect to $f$.

Then there exists a bound,
\begin{equation}
\Psi(\omega)=\Psi(\omega | n, f, a, b, \alpha, \alpha', \xi,
 \Lambda, diam(r^{-1}(a,b))
\end{equation}
such that
$$
\lim_{\omega \to 0} \Psi(\omega | n, f, a, b, \alpha, \alpha', \xi,
 \Lambda, diam(r^{-1}(a,b))=0.     
$$ 
and there exists a length space, $X$, such that 
$$
d_{GH}(r^{-1}((a+\alpha, b-\alpha)),X\times_f(a+\alpha, b-\alpha))
\le \Psi(\omega)
$$
where the region $r^{-1}((a,b))$ is endowed with the localized distance 
functions, $d^{\alpha, \alpha'}$. 
\end{theorem}

\begin{note} \label{DefnX}
The length space, $X$, defined by Cheeger and Colding
is a length space defined to be arbitrarily close to the
set, $r^{-1}(a + \alpha')$, endowed with a localized
distance function $d_U$ with 
\begin{equation} 
U=r^{-1}(a+\alpha'-\xi,a+\alpha'-\xi). 
\end{equation}
See [ChCo, Prop 3.3].
\end{note}

\begin{note} \label{DefnFGH} 
More important for our purposes is the fact that the
Gromov-Hausdorff map for this theorem is defined,
\begin{equation} 
F_{GH}(x)=(\pi(x), r(x)).
\end{equation}
Here $\pi(x)=f_X ( \bar{\pi}(x))$,
where $f_X: r^{-1}(a + \alpha') \mapsto X$ is a Gromov-Hausdorff
map and $\bar{\pi}(x)$ is any point in $r^{-1}(\alpha+\alpha')$
closest to $x$ [ChCo Thm 3.6].   Note that $F_{GH}$ is not a uniquely
determined function, nor is in continuous.
\end{note}

\begin{note} \label{GHOnto}
In the process of proving that $F_{GH}$ is almost onto, see Lemma 3.38 of 
[ChCo], Cheeger and Colding prove a formula which implies that for any
$t \in (a +\alpha, b-\alpha)$, the restricted function,
\begin{equation} \label{LevelGH} 
F_{GH}: r^{-1}(t) \longmapsto X' \times_f \{t\}
\end{equation}
is almost onto.  (See [So1, Section 4.2] for more details).
Thus $X'$ and $X$ can be best described as being Gromov-Hausdorff close
to any given level set $r^{-1}(t)$ rescaled by $f(t)$ with the localized
distance function $d_U$ of Theorem~\ref{ChCo4.85}. 
\end{note}

The fact that $F_{GH}$ is a Gromov-Hausdorff map between level sets
will be crucial to our proof of Theorem~\ref{SublinDiam}. 

\begin{note} \label{PsiK}
The estimating function 
$\Psi(\omega | n, f, a, b, \alpha, \alpha', \xi,
\Lambda, diam(r^{-1}(a,b))$ of the Cheeger-Colding Theorem 
depends on the warping function, $f$, only through the following quantities: 
\begin{eqnarray}
K_1\,\,\ge\,\,\sup_{r\in [a,b]}\,\, |f(r)| \qquad
 & 
K_3\,\,\ge\,\sup_{r\in [a,b]}\,|f'(r)| 
 \\
K_2\,\ge\,\sup_{r\in [a,b]} \,\Bigg|\frac 1 {f(r)} \Bigg|  \qquad
 &
K_4\,\ge\,\sup_{r\in [a,b]}\, \Bigg|\frac {f''(r)} {f(r)} \Bigg|
\end{eqnarray}
Note that we normalize $f(r)$ so that $f(a)=1$.  See [So1, Sections 4.1-4.2]
for details.
\end{note}

\begin{remark} \label{calV} {\em
Cheeger and Colding do not state this theorem exactly as we have written 
it above.  In their statement, the functions $\Psi$, $N$ and $D$ do not depend 
on $diam(r^{-1}(a,b))$ but instead on a function $\cal{V}$.     
\begin{equation}
{\cal{V}}(u)=\inf_{q\in r^{-1}(a,b)} \frac{Vol(B_u(q))}{Vol(r^{-1}(a,b))}
\end{equation}
By examining Prop. 4.50 and Lemmas 3.28 and 3.32 of [ChCo], where
the dependence on $\cal{V}$ is introduced, 
and Prop 2.24 of [ChCo], which describes the
properties of this function, it is clear that this dependence
can be replaced by dependence on the minimum Ricci curvature, $\Lambda$,
the dimension, $n$, and the diameter of the set measured with respect to
the standard metric, $diam(r^{-1}(a+\alpha,b-\alpha)$.
The restatement in Theorem~\ref{ChCo4.85} is convenient for our
purposes.  }
\end{remark}


\subsection{The Asymptotic Almost Rigidity of Manifolds with 
Minimal Volume Growth}
\Pind

We can now state our asymptotic almost rigidity theorem.
Recall the Definition~\ref{MinVolGrowth} of minimal volume growth.
Recall the constant $V_\infty$ of (\ref{Vinfty}).

\begin{theorem}  \label{GHWarp}
Let $M^n$ be a manifold with a ray $\gamma$, nonnegative Ricci
curvature everywhere, 
\begin{equation}
Ricci(x) \ge \frac {(n-1)(\frac 1 4 - v^2)} {b(x)^2}
\textrm { on } b^{-1}([r_0, \infty))
\end{equation}
where $v\in(0,1/2]$ and minimal volume growth.

Then for any given $\varepsilon>0$ and $L>\vare>0$, there exists 
a sufficiantly large constant, $V_{\vare,L} <V_\infty$,
such that if
\begin{equation}
Vol_{(n-1)}(b^{-1}(r_1)) \ge V_{\varepsilon,L} (r_1)^{(n-1)(1/2 -v)}
\end{equation}
then there exists a length space $X_{r_1}$ such that 
$$
d_{GH} \bigg(b^{-1}\Big((r_1+\varepsilon,r_1+L)\Big),
   X_{r_1+\vare} \times_{(b^{(1/2-v)})} (r_1+\varepsilon,r_1+L)\bigg) 
< \varepsilon \, diam\big(b^{-1}(r_1)\big).
$$
This Gromov-Hausdorff closeness is from level set to level set
[Note~\ref{GHOnto}], so in fact
$X_{r_1}$ is Gromov-Hausdorff close to any level $b^{-1}(s)$ rescaled
by the warping function $f(s)$ as long as $s\in (r_1+\vare,r_1+L)$. 

The distance function on $b^{-1}((r_1+\varepsilon,r_1+L))$
and $b^{-1}(s)$ is the localized distance function, 
$d_{b^{-1}(r_1+\varepsilon/2, r_1+L+\varepsilon/2)}$. 

\end{theorem}

The following is a corollary of the above or can be proven
directly with a simplification of the above theorem's proof.

\begin{theorem}  \label{GHCross}  Given a manifold, $M^n$,
with nonnegative Ricci curvature and linear volume growth,
for any given $\varepsilon>0$ and $L>\vare>0$, there exists 
a sufficiantly large constant,
\begin{equation}
V_{\vare,L} < \lim_{R\to\infty} Vol_{n-1}(b^{-1}(R)) =V_\infty
\end{equation}
such that if
\begin{equation}
Vol_{(n-1)}(b^{-1}(r_1)) \ge V_{\varepsilon,L}>0
\end{equation}
then 
$$
d_{GH} \bigg(b^{-1}\Big((r_1+\varepsilon,r_1+L)\Big),
   X_{r_1} \times 
(r_1+\varepsilon,r_1+L)\bigg) 
< \varepsilon \, diam\big(b^{-1}(r_1)\big).
$$
Here $X_{r_1}$ is a length space such that
$
d_{GH}\left(X_{r_1},b^{-1}\left(r_1+ \varepsilon \right)\right)<\varepsilon
\, diam\big(b^{-1}(r_1)\big).
$ 
All spaces in this theorem are endowed with the localized distance function
$d_U$ with $U=b^{-1}(r_1+\varepsilon/2, r_1+L+\varepsilon/2)$. 
\end{theorem}

Theorem~\ref{GHCross} 
essentially asserts that once a level set has a large enough 
$(n-1)$-volume, then the nearby region is almost an isometric
product of that level with an interval.  Note that
we are forced to shift our region over slightly
in order to be able to match the components of the region to that of
the level.  Cheeger and Colding are only able to control the distances
of a subregion of the original region because the estimates on the
Hessian of the distance function are only controlled
on subregions of the region where the volume is controlled.  
For this reason, all the 
distance functions are also localized inside subsets of the original region.

Note that these manifolds do not necessarily converge to unique
warped product manifolds even if the diameter of the Busemann level
sets is uniformly bounded.  In [So2], there are examples of manifolds
satisfying the hypothesis of Theorem~\ref{GHWarp} for which
there exist $r_i \to \infty$ such that 
\begin{equation}
b^{-1}((r_{2i}, r_{2i} +L)) \to X \times (0,L),
\end{equation}
and
\begin{equation}
b^{-1}((r_{2i+1}, r_{2i+1} +L)) \to Y \times (0,L),
\end{equation}
where $X$ and $Y$ are not isometric.
In order to force the manifold to be asymptotically close to a unique
isometric product manifold we would have to add additional conditions on 
the speed at which $Vol_{n-1}(b^{-1}(r))$ approaches $V_\infty$.
See Remarks~\ref{MoreWarp1} and~\ref{MoreWarp2} after the proof
of Theorem~\ref{GHWarp}.

\subsection{Minimal Volume Growth and Almost Maximality}

\Pind
In this section we begin a proof of both Theorem~\ref{GHWarp} on asymptotic
almost rigidity and Theorem~\ref{SublinDiam} on the diameter growth
of manifolds with minimal volume growth and globally nonnegative Ricci 
curvature.  The key ingredient in both proofs is the application of 
the Cheeger Colding Theorem~\ref{ChCo4.85}.  Here we provide a series of
lemmas which relate our hypothesis of
{\em minimal volume growth} to their hypothesis of {\em almost
maximality}. Recall Definitions~\ref{MinVolGrowth} and~\ref{AlmostMax}.

In order to study noncompact manifolds with minimal volume growth, we do
not examine standard distance functions, but instead we examine
Busemann functions.   The following lemma asserts that Busemann functions
are distance functions on certain regions (See [So1] for the proof).

\begin{lemma} \label{BuseDist}  
Let $r_2$ be a real number.  Let $r(x)=d(x, b^{-1}(r_2))$.
Then $r(x)=r_2-b(x)$ on the region $b^{-1}(-\infty, r_2]$.
\end{lemma}

Recall that Theorem~\ref{CompactLevels} states that 
the Busemann level sets are compact on manifolds with
minimal volume growth and globally nonnegative Ricci curvature.  
Thus $r(x)=d(x, b^{-1}(r_2))$ of the above lemma can be used
as our distance function in the Cheeger-Colding Theorem [ChCo, 4.85].

Cheeger and Colding showed that regions with  
almost maximal volume were almost warped products.  Here, we are
studying compact regions in noncompact manifolds with minimal volume growth.
These ideas are related because
in a manifold with minimal volume growth, annuli about increasingly
distant points, $\gamma(R_i)$, have almost maximal volume.  
Such annuli converge to regions between level sets of the Busemann function,
$b_\gamma$.  See Lemma~\ref{BuseInnerAnn} and [So2]. 

In the next lemma we show that once a Busemann level set in such 
noncompact manifold has sufficiently large (n-1)-volume, then any region beyond
that level is $\omega$-almost maximal.
Recall, also the definition of the $\omega$-almost
maximal volume property in Definition~\ref{AlmostMax}
and of $V_\infty$ in (\ref{Vinfty}).

\begin{lemma}  \label{MinVolToAlmostMax}  
Let $M^n$ be a manifold with a ray $\gamma$ and Ricci curvature bounded below
as in (\ref{RicciBound}) with $v\in [0, (n+1)/2(n-1))$.
Suppose $M^n$ has minimal volume growth and that
the level sets of the Busemann function are compact.

If $r_1$ is large enough that
\begin{equation}
0< V_\infty - \frac{Vol_{n-1}(b^{-1}(r_1))}{r_1^p} < \omega V_\infty
\end{equation}
then for all $r_2>r_1$, the region,
$
b^{-1}((r_1, r_2))=r^{-1}(0,r_2-r_1),
$ 
where $r(x)=r_2-b(x)=d(x,b^{-1}(r_2))$, has the $\omega$-almost
maximal volume property with respect to the function 
$f(r)=(r_2-r)^{(\frac 1 2 -v)}$.
\end{lemma}

\noindent
{\bf Proof:} 
First, it is easy to check that $f(r)$ is the appropriate function
for the Ricci bound because 
$$  
-(n-1) \frac {f''(r(x))}{f(r(x))}= -(n-1)
\left(\frac 1 2 -v\right)\left(-\frac 1 2 -v\right)
\frac {(-1)^2} {(r_2-r)^2}= (n-1)\frac{(\frac 1 4 -v^2)}{b^2(x)}.
$$
As for the mean curvature bound we know 
\begin{equation}
(n-1) \frac {f'(0)}{f(0)}= (n-1)\frac {-(\frac 1 2 -v)(r_2)^{(-\frac 1 2 -v)}}
     {(r_2)^{(\frac 1 2 -v)}}=\frac{-(n-1)(\frac 1 2 -v)}{r_2}
\end{equation}
Using the Laplacian Comparison Theorem [Ch]
and (\ref{eqnlapb}) which tells us the Laplacian of the Busemann
function on our comparison warped product manifold, we get,
\begin{equation}
\frac{-(n-1)(\frac 1 2 -v)}{r_2}
     \ge-\Delta b = -m( b^{-1}(r_2))
     =m(r^{-1}(0)).
\end{equation}
Thus we have verified the mean curvature requirement.

So now we need only show 
\begin{equation}
\frac{Vol(r^{-1}(0,b))}{Vol_{n-1}(r^{-1}(0))}\,\,\ge 
\,\,\,(1-\omega)\,\frac{\int_0^b \, f^{n-1}(r)\,\,dr}{f^{n-1}(0)},  
\end{equation}
or equivalently,
\begin{equation}
\frac{Vol(b^{-1}(r_1,r_2))}{Vol_{n-1}(b^{-1}(r_2))}\,\,\ge\,\, 
(1-\omega)\,\frac{\int_{r_1}^{r_2} \, b^{(n-1)(\frac 1 2 -v)}\,\,db}
{(r_2)^{(n-1)(\frac 1 2 -v)}}.
\end{equation}
By Theorem~\ref{Monotonicity} we have
\begin{eqnarray*}
\frac{\,Vol(b^{-1}(r_1,r_2))\,}{\,Vol_{n-1}(b^{-1}(r_2))\,} & = &
         \frac{\,\int_{r_1}^{r_2}\, Vol_{n-1}(b^{-1}(s))\,\,ds\,}
                      {Vol_{n-1}(b^{-1}(r_2))} \\
 &\ge& \frac{\,\bigg(\int_{r_1}^{r_2}\, 
       \big(\frac s {r_1}\big)^{(n-1)(\frac 1 2 -v)}\,\,ds\bigg)
         \,\, Vol_{n-1}\big(b^{-1}(r_1)\big)\,}
         {Vol_{n-1}\big(b^{-1}(r_2)\big)}\\
 & = & \frac {\,Vol_{n-1}(b^{-1}(r_1))\,}{Vol_{n-1}(b^{-1}(r_2))}
      \, \frac{\,\int_{r_1}^{r_2} \, b^{(n-1)(\frac 1 2 -v)}\,\,db\,}
             {(r_2)^{(n-1)(\frac 1 2 -v)}}
       \bigg(\frac{r_2}{r_1}\bigg)^{(n-1)(\frac 1 2 -v)}.
\end{eqnarray*}
So we need only show 
\begin{equation}  \label{ShowMon1}
\frac {\,\,Vol_{n-1}(b^{-1}(r_1))\,\,}{Vol_{n-1}(b^{-1}(r_2))}\,
        \bigg(\frac{r_2}{r_1}\bigg)^{(n-1)(\frac 1 2 -v)} \,\ge \,(1-\omega).
\end{equation}
Now we are given 
\begin{equation}
\frac {\,\,\, \left(V_\infty-
     \frac{Vol_{n-1}(b^{-1}(r_1))}{r_1^{(n-1)(\frac 1 2 -v)}}\right)\,\,\, }
       { V_\infty }\,\, \le\,\, \omega
\end{equation}
Using the monotonicity,
\begin{equation}
\frac{\,Vol_{n-1}(b^{-1}(r_1))\,}{r_1^{(n-1)(\frac 1 2 -v)}}\,\,<\,\,
\frac{\,Vol_{n-1}(b^{-1}(r_2))\,}{r_2^{(n-1)(\frac 1 2 -v)}}
\,\,\le \,\,V_\infty,
\end{equation}
we find
\begin{equation}
\frac{\,\,\left( \frac{Vol_{n-1}(b^{-1}(r_2))}{r_2^{(n-1)(\frac 1 2 -v)}}-
       \frac{Vol_{n-1}(b^{-1}(r_1))}{r_1^{(n-1)(\frac 1 2 -v)}}\right) \,\,}
  {\left( \frac{Vol_{n-1}(b^{-1}(r_2))}{r_2^{(n-1)(\frac 1 2 -v)}}\right)} 
  \,\,\,\le \,\,\,\omega.
\end{equation}
Cancelling the terms involving $r_2$ and rearranging
this equation, we obtain (\ref{ShowMon1}) and we are done.

\ProofEnd

In the Cheeger Colding Theorem~\ref{ChCo4.85}, 
the given region has a fixed lower Ricci curvature bound, diameter
bound and comparison warping function.
It is shown that the region is almost a warped
product if its almost maximal volume estimate, $\omega$, is sufficiently close
to $0$.   The Gromov-Hausdorff closeness depends on 
\begin{equation}
\Psi(\omega_{\vare, L} | n, f, a, b, \alpha, \alpha', 
\xi,\Lambda, diam(r^{-1}(a,b)))
\end{equation} 
and $\Psi$ only approaches $0$ when all the other parameters are fixed.
 
In our situation, we are
examining the asymptotic behavior of a sequence of regions contained between 
Busemann levels, $b^{-1}(r_1,r_2)$, where $r_1$ and $r_2$ approach infinity.
Thus our set $K=b^{-1}(r_2)$ will not
be a fixed set and the diameters of the regions, $b^{-1}(r_1,r_2)$, will be
changing.  Thus we must rescale the regions before applying
the Cheeger-Colding Theorem.  

\begin{lemma} \label{Rescale} 
Let $N^n_1$ be a Riemannian manifold with a region $r_1^{-1}(a_1,b_1)$ that has
$\omega$-almost maximal volume with respect to a
function $f_1$.
Let $N^n_2$ be the manifold $N^n_1$ with its metric scaled down by $D^2$.
Then $r_2(x)= (r_1(x)-a_1)/D$ is a distance function on the region
$r_2^{-1}(0,(b_1-a_1)/D)$ in $N_2^n$.  Furthermore, if we let
$f_2(t)=f_1(tD+a_1)/D$ then the region, 
$r_2^{-1}((0, \frac{b_1-a_1}{D}))$, has
$\omega$-almost maximal volume with respect to this function $f_2$.
\end{lemma}

The details of this proof can be found in [So1, Lemma 4.10]. 

Clearly, the rescaling of the manifold will affect other parameters in the
Cheeger-Colding theorem.  In particular, the distance between levels sets,
$r^{-1}(a)$ and $r^{-1}(b)$ may become very small.  This is a problem
because the Cheeger-Colding Theorem requires that there be fixed constants
$b-a>\alpha>\alpha'>0$.  While we cannot employ the Cheeger-Colding Theorem
to prove the Gromov-Hausdorff closeness in this situation it is easy to see
that a thin set must be close to a warped product regardless of its volume 
properties.  See [So1, Lemma 4.11] for details.

\begin{lemma} \label{ThinGH}
Given a manifold, $M^n$, a compact subset, $K$,
and a distance function $r(x)=d(x,K)$.
Given any $\vare>0$, if
\begin{equation}
b-a<\delta_{\vare}=\frac{\vare}{4}
\end{equation}
then
\begin{equation}
d_{GH}(r^{-1}(a,b), r^{-1}(a)\times_f(a,b))<\vare.
\end{equation}
where the distance function on $r^{-1}(a,b)$ and $r^{-1}(a)$
can be any localized distance function,
$d_U$, where $U\supset r^{-1}(a,b)$.  (See Defn~\ref{DefnDist}).
\end{lemma}

Before rescaling regions between level sets, we would like to estimate their 
diameter.  Three different estimates are obtained in the following lemmas.
The first lemma is simple but is used to prove both Theorem~\ref{GHWarp}
and Theorem~\ref{SublinDiam}.

\begin{lemma} \label{RayDiam}
Let $M^n$ be any complete noncompact Riemannian manifold with
a Busemann function, $b$.  Then, for all $r_1<r_2$,
\begin{equation}
diam(b^{-1}(r_1, r_2)) \le diam(b^{-1}(r_2)) + 2(r_2-r_1).
\end{equation}
\end{lemma}

\noindent {\bf Proof:}
Given any $x$ in $b^{-1}(r_1)$, there exists a Busemann ray, $\gamma_x$,
which is parametrized by arclength,
such that $\gamma_x(r_1)=x$ and $\gamma_x(r_2) \in b^{-1}(r_2)$.
\ProofEnd

Thus to control the diameter of a region we need only control the
diameter of the boundary closer to infinity.

The next lemma gives an explicit bound on the diameter of the
boundary as a function of the diameter and volume of the first
level set.  This lemma cannot employ the Busemann rays to travel 
between the levels and thus requires minimal volume growth and 
globally nonnegative Ricci curvature.
The techniques used to prove this lemma were developed in [So2].

\begin{lemma} \label{DiamChange}  
Let $M^n$ be a manifold with $Ricci \ge 0$
everywhere.  Suppose it has  Ricci Curvature bounded as in (\ref{RicciBound})
with $v\in [0,1/4]$ and minimal volume growth.

Given any $0<\vare<1/2$, we can find a level set
with a sufficiently large volume
\begin{equation} \label{DiamChange1}
\frac{Vol_{n-1}(b^{-1}(r_1))}{pr_1^{p-1}} > V_\infty 
\bigg(1-\frac 1 2 \bigg(\frac{\vare}{2+2\vare}\bigg)^n\bigg)
\end{equation}
such for any $L>0$ and any $r_3>r_1$ the level sets $b^{-1}(r_3)$ and
$b^{-1}(r_3 + L)$ are Hausdorff close as subsets of $M^n$
$$
d_H\big(b^{-1}(r_3),b^{-1}(r_3 + L)\big)<
                  \vare\, diam\big(b^{-1}(r_3 + L)\big) + L
$$
and the difference between their diameters is controlled
\begin{equation}
|diam(b^{-1}(r_3))-diam(b^{-1}(r_3+L))|
<\vare \,diam(b^{-1}(r_3+L)) + 2L.
\end{equation}
\end{lemma}

\noindent
{\bf Proof:}  
By the last lemma, we know that
\begin{equation}  \label{DiamChange3}
b^{-1}(r_3) \subset T_\delta\big(b^{-1}(r_3+L)\big).
\end{equation}
Thus the diameter of the first level set can be no larger than that of the
second level set plus $2L$,
\begin{equation}
diam(b^{-1}(r_3))<diam(b^{-1}(r_3+L)) +2L.
\end{equation}

So now we must bound $diam(b^{-1}(r_3+L))$ from above and
show that it is contained in the appropriate tubular neighborhood
of $b^{-1}(r_3)$.
We will do this proof by contradiction.

Let $D=diam(b^{-1}(r_3 + L))$.  Suppose that 
\begin{equation} \label{DiamChange2}
b^{-1}(r_3+L) \notin
T_{(\vare D/2)}\big(\Omega_{r_3, r_3+L+\vare D/2}(b^{-1}(r_3))\big).
\end{equation}
Then there is a point $x\in b^{-1}(r_3+L)$ such that 
\begin{equation}
B_x(\vare D/2) \cap 
\Omega_{r_3, r_3+L+\vare D/2}(b^{-1}(r_3)))=0.
\end{equation}
Thus, by Theorem~\ref{Monotonicity}, we have
\begin{eqnarray*}
\lefteqn{Vol(B_x(\vare D/2))<}\\
& < & Vol(b^{-1}(r_3+L-\vare D/2,r_3+L+\vare D/2))\\
&   & \qquad
  -\qquad Vol(\Omega_{r_3+L-\vare D/2, r_3+L+\vare D/2}(b^{-1}(r_3))) \\
& < & \bigg( \bigg(r_3+L+\frac {\vare D}{R} \bigg)^p-
    \bigg(r_3+L-\frac{\vare D}{R}\bigg)^p \bigg)
          \bigg(V_\infty
              -  \frac {Vol(b^{-1}(r_3))}{pr_3^{p-1}}\bigg).
\end{eqnarray*}

On the other hand, by the Relative Volume Comparison Theorem [Bi][GrLaPa]
and $Ricci\ge 0$ everywhere, we know
\begin{equation} \label{huge}
Vol(B_x(\vare D/2))>Vol(B_x(R))\bigg(\frac {\vare D/2}{R}\bigg)^{n}
\end{equation}
We set $R= diam(b^{-1}(r_3+L))+ 2 (\vare D/2)= D+2 (\vare D/2)$ to insure that
the ball of radius
$R$ contains $b^{-1}(r_3+L-\vare D/2,r_3+L)$.
Thus
$$ 
Vol(B_x(\vare D/2)) \,\, >\,\,  Vol\big(b^{-1}(r_3+L-\vare D/2,r_3+L)\big)
                           \bigg(\frac {\vare D/2} R \bigg)^{n}  
$$
$$
  \qquad  > \,\,\frac{ Vol_{n-1}\big(b^{-1}(r_3)\big)}{pr_3^{p-1}}
   \bigg( (r_3+L+\frac {\vare D}{R} )^p-(r_3+L-\frac{\vare D}{R})^p \bigg)
                        \bigg(\frac {\vare D/2}{D+2 (\vare D/2)} \bigg)^{n}.
$$
This last line employs Theorem~\ref{Monotonicity} once again.

Using these two bounds for $Vol(B_x(\vare D/2))$, we have
$$ 
\bigg( \bigg(r_3+L+\frac {\vare D}{R} \bigg)^p-
   \bigg(r_3+L-\frac{\vare D}{R}\bigg)^p \bigg) 
   \left(V_\infty-\frac{Vol_{n-1}(b^{-1}(r_3))}{pr_3^{p-1}}\right)\,\,\,  > 
  \qquad
$$
$$
\qquad >\,\,\, \frac{Vol_{n-1}(b^{-1}(r_3))}{pr_3^{p-1}}
\bigg( (r_3+L+\frac {\vare D}{R} )^p-(r_3+L-\frac{\vare D}{R})^p \bigg)
             \bigg(\frac{\vare D/2}{D(1+\vare)}\bigg)^{n}
$$
which we can rewrite as
\begin{equation}
\frac{(V_\infty p r_3^{p-1}-Vol_{n-1}(b^{-1}(r_3)))}{Vol_{n-1}(b^{-1}(r_3))}
        >  
     \bigg(\frac {\vare/2}{(1+\vare)}\bigg)^{n}\frac 1 2 . 
\end{equation}
Now we take $r_3$ large enough for our volume estimate given by our choice
of $r_1$ in (\ref{DiamChange1}) and get
\begin{eqnarray}
\frac 1 2 \bigg(\frac{\vare}{2+2\vare}\bigg)^n &>& 
\frac{V_\infty p r_3^{p-1} -Vol_{n-1}(b^{-1}(r_1))}{V_\infty}\\
&> &\frac{(V_\infty p r_3^{p-1}-Vol_{n-1}(b^{-1}(r_3)))}
   {Vol_{n-1}(b^{-1}(r_3))}\\
&>& \bigg(\frac {\vare/2}{(1+\vare)}\bigg)^{n}\frac 1 2 
\end{eqnarray}
which is a contradiction.

Thus our assumption in (\ref{DiamChange2}) does not hold and
instead we have
\begin{equation} \label{DiamChange4}
b^{-1}(r_3+L) \subset  
T_{(\vare D/2)}\bigg(\Omega_{r_3, r_3+L+\vare D/2}(b^{-1}(r_3))\bigg).
\end{equation}
Since all points in $\Omega_{r_3, r_3+L+\vare D/2}(b^{-1}(r_3))$ 
are on segments of Busemann rays running from $b^{-1}(r_3)$ 
of length less than or equal to $L+\vare D/2$, we have
\begin{equation}\label{DiamChange15}
\Omega_{r_3, r_3+L+\vare D/2}(b^{-1}(r_3)))\subset 
T_{(L+\vare D/2)}\bigg(b^{-1}(r_3)\bigg).
\end{equation}
So, combining (\ref{DiamChange4}) and (\ref{DiamChange15}) we get
\begin{equation}
b^{-1}(r_3+L) \subset T_{(L+\vare D/2+\vare D/2)}\big(b^{-1}(r_3)\big),
\end{equation}
where $D=diam(b^{-1}(r_3 + L))$.

The lemma then follows.
%

\ProofEnd

We end this section with an easy rough estimate for a lower bound 
on the diameter of a Busemann level set as a function of its
(n-1)-volume.  Once again we restrict ourselves to manifolds with
nonnegative Ricci curvature everywhere.
However, we do not assume that we have minimal volume growth.

\begin{lemma} \label{LowBoundDiam}
In a complete noncompact manifold with nonnegative Ricci curvature
such that $b^{-1}(r)$ is compact, we have
\begin{equation} 
diam(b^{-1}(r)) > \left(\frac{Vol_{n-1}(b^{-1}(r))}
                    {w_n\,2^n}\right)^{(1/(n-1))}.
\end{equation}
\end{lemma}

\noindent
{\bf Proof:}
Let $d=diam(b^{-1}(r))$.

For any $x\in b^{-1}(r)$, $B_x(2d) \supset b^{-1}(r,r+d)$, so
\begin{equation}
Vol(B_x(2d)) > Vol(b^{-1}(r,r+d))
\end{equation}
By the Bishop Volume Comparison Theorem [Bi, BiCr],
\begin{equation}
w_n(2d)^n \ge Vol(B_x(2d)).
\end{equation}
Since $b$ is Lipschitz we can employ the Coarea formula [Fed 3.2.11] 
and Lemma~\ref{Monotonicity} to get,
\begin{equation}
Vol(b^{-1}(r,r+d)) \ge (d) Vol_{n-1}(b^{-1}(r)).
\end{equation}
Thus
\begin{equation}
w_n(2d)^n > (d) Vol_{n-1}(b^{-1}(r))
\end{equation}
and
\begin{equation}
d^{n-1} >\frac{Vol_{n-1}(b^{-1}(r))}{w_n 2^n}.
\end{equation}
\nopagebreak[4]
\ProofEnd

\subsection{The Proof of the Asymptotic Almost Rigidity Theorem}

We now prove Theorem~\ref{GHWarp}.  Throughout
this section  $M^n$ satisfies the hypotheses of
this theorem.  See Section 2.2.   

We will begin by rescaling the region,
$b^{-1}(r_1, r_1+L+ \vare)$, so that the diameter of the rescaled region 
is bounded above.  Then we can apply the Cheeger Colding Theorem.
Since, 
\begin{equation}
diam(b^{-1}(r_1, r_1+L +\vare) \le 2(L +\vare) + diam(b^{-1}(r_1+\vare)
\end{equation}
by Lemma~\ref{RayDiam}, we divide the metric by $D^2$ where 
\begin{equation}\label{GenWarpD}
D=diam(b^{-1}(r_1 +L +\vare)).
\end{equation}

We will first prove that given any $\vare >0$ and any $L>\vare>0$ there exists
$V_{\vare, L}$ sufficiently large that if
\begin{equation}
Vol_{(n-1)}(b^{-1}(r_1)) \ge V_{\varepsilon,L} (r_1)^{(n-1)(1/2-v)}
\end{equation}
then 
\begin{equation} \label{NotFinal}
d_{GH} \bigg(b^{-1}\Big((r_1+\varepsilon,r_1+L)\Big),
   X_{r_1+\vare} \times_{(b^{(1/2-v)})} (r_1+\varepsilon,r_1+L)\bigg) 
< \varepsilon' \, D,
\end{equation}
where
\begin{equation} \label{ReNormEps}
\vare ' \,\,\le\,\, 
\frac{\vare\,\,\,\left(\frac {V_\infty}{2^{n+1}}\right)^{1/(n-1)}} 
{\,\,2 \left(\frac {V_\infty}{2^{n+1}}\right)^{1/(n-1)}
+\,4L\,+\,4\vare\,\,}.
\end{equation}
This strange choice of $\vare'$ has been made so that later we can replace
the dependence on $diam\big(b^{-1}(r_1+L+\vare)\big)$ by
$diam\big(b^{-1}(r_1)\big)$ using Lemma~\ref{DiamChange}.

We let $r(x)=d(x, b^{-1}(r_1 +L + \vare))=(r_1+L+ \vare -b(x))/D$
in this rescaled region.  Note that $r(x)$ increases as $b(x)$ decreases
and that this region is thus $r^{-1}(0,b)$ where 
\begin{equation}\label{GenWarpB}
b=\frac{L +\vare}{D}.
\end{equation}

When $b< 2\vare'$, the region is ``thin'' and we know
\begin{equation}\label{LessFinal}
d_{GH}(r^{-1}(a,b), r^{-1}(a)\times_f(a,b))<\vare'.
\end{equation}
for any warping function $f$ by Lemma~\ref{ThinGH}.
Thus (\ref{NotFinal}) holds when we rescale back to the original region.

So we will assume $b\ge 2\vare$ and apply the Cheeger-Colding Theorem
to obtain (\ref{LessFinal}).  First we must bound all the parameters in 
\begin{equation}
\Psi(\omega | n, f, a, b, \alpha, \alpha', \xi,\Lambda, diam(r^{-1}(a,b))),
\end{equation} 
of the Cheeger-Colding theorem for our rescaled region.
We set $n=n$, the dimension of our manifold.  Since 
$Ricci \ge 0$ globally, we have $\Lambda =0$.  We have $a=0$ and 
\begin{equation}
b=\frac{L +\vare}{D} \le 
\frac{L+\vare}
   {\,\,\,\,\bigg(\frac{Vol_{n-1}(b^{-1}(r_1+L+\vare)}{2^n}\bigg)^{(1/(n-1)}}
\end{equation}
by Lemma~\ref{LowBoundDiam} which provides a lower bound on the diameter
in terms of the volume.
If we take 
\begin{equation}
Vol_{n-1}(b^{-1}(r_1 + L +\vare))\ge V_{\varepsilon,L}>V_\infty/2,
\end{equation}
we can bound $b$ from above by a constant 
\begin{equation}
b\le\frac{L+\vare}{\bigg(\frac{V_\infty/2}{2^n}\bigg)^{(1/(n-1))}}.
\end{equation}
The other parameters can be set $\alpha'=\alpha/2$ and $\chi=\alpha'/2$
where
\begin{equation}
\alpha=\frac{\vare}{diam(b^{-1}(r_1 +L +\vare))}
> \frac {b\,\vare} {(L+\vare)}.
\end{equation}
so they are all bounded from below and above as well.

The warping function, $f$, normalized such that $f(0)=1$, is
\begin{equation}
f(r)=\left( \frac{r_1+Db-Dr}{r_1+Db} \right)^{\left(\frac 1 2 - v\right)}.
\end{equation}
To employ the Cheeger-Colding Theorem, we must bound this warping function
as described in Remark~\ref{PsiK}.  More specifically, 
to bound the constants, $K_i$,
uniformly for all values of $r_1$ since we will have to vary $r_1$
to obtain the $\omega$-almost maximality required by the Cheeger-Colding 
Theorem.  This is especially troublesome because $D=diam(b^{-1}(r_1+L+\vare))$
depends on $r_1$ and may approach infinity.

To bound the constants $K_i$ which depend on $f$, 
we will use the fact that the diameter of the Busemann
levels grows at most linearly [Theorem~\ref{LinDiam}].
We must chose
\begin{equation}
K_1 =  \sup_{r\in [0,b]} |f(r)|= \left( \frac{r_1+Db}{r_1+Db} \right)^{\left(\frac 1 2 - v\right)}=1.
\end{equation}
We need
\begin{equation}
K_2 \ge \sup_{r\in [0,b]} \Bigg|\frac 1 {f(r)} \Bigg| 
=\left( \frac{r_1}{r_1+Db} \right)^{-\left(\frac 1 2 - v\right)}.
\end{equation}
Now $b$ is bounded above and below so we need only worry about $r_1$
and $D$ as $r_1$ approaches infinity.  By the at most linear diameter growth,
we have $D_{r_1}\le C r_1$; so
\begin{equation}
\lim_{r_1\to\infty}
\left( \frac{r_1}{r_1+Db} \right)^{-\left(\frac 1 2 - v\right)}
\le \left(\frac 1 {1+Cb}\right)^{-\left(\frac 1 2 - v\right)}.
\end{equation} 
and $K_2$ exists.
We must chose $K_3$ such that
\begin{eqnarray*}
K_3 &\ge&\sup_{r\in [0,b]}|f'(r)|\\
&=&\sup_{r\in [0,b]}\bigg|{\left(\frac 1 2 - v\right)}
\left( \frac{r_1+Db-Dr}{r_1+Db} \right)^{\left(-\frac 1 2 - v\right)}
\left(\frac{-D}{r_1+Db}\right)\bigg|\\
&=&{\left(\frac 1 2 - v\right)}\left(\frac{D}{r_1+Db}\right)
\left( \frac{r_1}{r_1+Db} \right)^{\left(-\frac 1 2 - v\right)}.
\end{eqnarray*}
Once again we verify that the right hand side is bounded as $r_1$ 
goes to infinity,
\begin{eqnarray*}
\lefteqn{\lim_{r_1\to\infty}{\left(\frac 1 2 - v\right)}
\left(\frac{D}{r_1+Db}\right)
\left( \frac{r_1}{r_1+Db} \right)^{\left(-\frac 1 2 - v\right)}\le}\qquad\\
\qquad&\le& {\left(\frac 1 2 - v\right)}\left(\frac{C}{1}\right)
\left( \frac{1}{1+Cb} \right)^{\left(-\frac 1 2 - v\right)}.
\end{eqnarray*}
Finally we need to chose $K_4$ such that
\begin{eqnarray*}
K_4 &\ge &\sup_{r\in [0,b]} \Bigg|\frac {f''(r)} {f(r)} \Bigg|\\
&=&\sup_{r\in [0,b]}{\left(\frac 1 2 - v\right)}{\Bigg|-\frac 1 2 - v\Bigg|}
     \left(\frac{D}{r_1+Db}\right)^2
\left( \frac{r_1+Db-Dr}{r_1+Db} \right)^{-2}\\
&=&{\left(\frac 1 2 - v\right)}{\left(\frac 1 2 + v\right)}
     \left(\frac{D}{r_1+Db}\right)^2
\left( \frac{r_1}{r_1+Db} \right)^{-2}.\\
\end{eqnarray*}
Once again we check if this is bounded as $r_1$ approaches infinity 
\begin{eqnarray*}
\lefteqn{
\lim_{r_1\to\infty}{\left(\frac 1 2 - v\right)}{\left(\frac 1 2 + v\right)}
     \left(\frac{D}{r_1+Db}\right)^2
\left( \frac{r_1}{r_1+Db} \right)^{-2}\le}\qquad\\
\qquad&\le& {\left(\frac 1 2 - v\right)}{\left(\frac 1 2 + v\right)}
\left(\frac{C}{1}\right)^2
\left( \frac{1}{1+Cb} \right)^{-2}
\end{eqnarray*}
Thus we have uniformly fixed our constants $K_i$.

For fixed $\vare$ and $L$ and keeping all our parameters bounded 
as above, we can find an $\omega_{\vare,L}$, 
depending only on
$\vare$ and $L$ such that
\begin{equation} \label{eqnGenPsi}
\Psi(\omega_{\vare, L} | n, f, a, b, \alpha, \alpha', \xi,\Lambda, 
diam(r^{-1}(a,b)))<\vare',
\end{equation}
where $\vare'$ is defined in (\ref{ReNormEps}).

If we can show that $r^{-1}(a,b)$ is $\omega_{\vare,L}$-almost maximal, then by
the Cheeger-Colding Theorem we have
\begin{equation} \label{gen1}
d_{GH}\big(r^{-1}(a+\alpha, b-\alpha), 
Y\times_f(a+\alpha, b-\alpha)\big)\, < \,\vare', 
\end{equation}
where $Y$ is a length space and each level $b^{-1}(s)$ 
is in fact mapped by the 
Gromov-Hausdorff equivalence map almost onto $f(s)Y$ [Note~\ref{GHOnto}].  
The distance
function on the metric product and the region are the localized distance
functions $d^{\alpha,\alpha'}=d_{r^{-1}(a+\alpha', b-\alpha')}$ defined 
in Definition~\ref{DefnDist}.

By the Lemma~\ref{Rescale} we need only show that $b^{-1}(r_1,R_1+L+\vare)$
has $\omega_{\vare,L}$-almost maximal volume to show that the rescaled region
$r^{-1}(a,b)$ is $\omega_{\vare,L}$-almost maximal as well.  By 
Lemma~\ref{MinVolToAlmostMax}, we know that if we take $r_1$ large enough
that $Vol_{n-1}(b^{-1}(r_1))\ge V_{\vare,L}$ where 
\begin{equation} \label{genV0}
V_{\vare, L}>(1-\omega_{\vare,L})V_\infty
\end{equation}
then $b^{-1}(r_1,R_1+L+\vare)$ has $\omega_{\vare,L}$-almost maximal volume.
Thus (\ref{gen1}) holds.  

If we rescale the region $r^{-1}(a+\alpha,b-\alpha)$ back up to
the original region $b^{-1}(r_1+\vare, r_1+L)$, then we can rescale
(\ref{gen1}) to get
\begin{equation} \label{NotFinal2}
d_{GH} \bigg(b^{-1}\Big((r_1+\varepsilon,r_1+L)\Big),
   X_{r_1+\vare} \times_{(b^{(1/2-v)})} (r_1+\varepsilon,r_1+L)\bigg) 
< \varepsilon' \, D,
\end{equation}
where the distance function on these spaces is the localized distance
function, $d_{b^{-1}(r_1+\vare/2, r_1+L+\vare/2)}$.  The distance function
rescales in  this manner because $\alpha'=\alpha/2=\vare/(2D)$.
Since this closeness holds on each level set, $X_{r_1+\vare}$
is close to $b^{-1}(r_1+\vare)$.

Thus we have obtained (\ref{NotFinal}).
To complete the proof of our theorem we need to show that (\ref{NotFinal2})
holds if we replace the $\vare' diam(b^{-1}(r_1+L+\vare))$ by 
$\vare diam(b^{-1}(r_1))$.

This may require us to take $r_1$ a little further out so that
the volume of its level set is close enough to $V_\infty$ to employ 
Lemma~\ref{DiamChange}.
That is, we take
\begin{equation} \label{eqnV3}
\frac{ Vol_{n-1}(b^{-1}(r_1)) }{ pr_1^{p-1} } > V_{\vare, L}
>V_\infty \bigg( 1-\bigg(\frac{1}{6}\bigg)^n\bigg)>\frac{V_\infty}{2},
\end{equation}
which implies that
\begin{equation}
V_\infty -\frac{Vol_{n-1}(b^{-1}(r_1))}{pr_1^{p-1}}<V_\infty
\bigg(\frac{1/2}{2+2(1/2)}\bigg)^n.
\end{equation}
So we by Lemma~\ref{DiamChange}, we know that
\begin{equation}
diam\big(b^{-1}(r_1+L+\vare)\big)\,<\,2\,diam(b^{-1}(r_1)) + 4(L+\vare).
\end{equation}

Using this information, we can rewrite our estimate
in (\ref{NotFinal}) as
\begin{eqnarray*}
\lefteqn{ d_{GH} \bigg(  {b^{-1}\bigg((r_1+\vare,r_1+L)\bigg)},
 { X_{r_1+\vare} \times_{b^{1/2-v}} \bigg(r_1+\vare,r_1+L\bigg)}
\bigg) < \,\,\,\,\,}  \\ 
\,\,& < &\,\,\frac{\vare \Dmin}{\bigg(2\Dmin + 4\bigg(L+\vare\bigg)\bigg)}
                      \,\,\,diam(b^{-1}(r_1+L+\vare))\\
\,\,& < &\frac{\vare \Dmin}{\bigg(2\Dmin + 4\bigg(L+\vare\bigg)\bigg)}  
                     \bigg(2diam(b^{-1}(r_1)) + 4(L+\vare)\bigg).
\end{eqnarray*}
By (\ref{eqnV3}), we have
$Vol_{n-1}(b^{-1}(r_1)) > V_\infty/2$; so by
Lemma~\ref{LowBoundDiam} we know that
\begin{equation}
diam\big(b^{-1}(r_1)\big) > \Dminr > \Dmin.
\end{equation}
Since $D/(2D+4(L+\vare))$ is a decreasing function of D, we can substitute
this diameter estimate in the Gromov-Hausdorff estimate to get
\begin{eqnarray*}
\lefteqn{ d_{GH} \bigg(  {b^{-1}\bigg(r_1+\vare,r_1+L\bigg)},
 { b^{-1}\bigg(r_1+\vare\bigg) \times \bigg(r_1+\vare,r_1+L\bigg)}
\bigg) < \,\,\,\,\,\,\,\,\,\,}  \\ 
\,\,\,\,\,\,& < &\frac{\vare diam(b^{-1}(r_1))}
         {\bigg(2diam(b^{-1}(r_1)) + 4\bigg(L+\vare\bigg)\bigg)}  
                     \bigg(2diam(b^{-1}(r_3)) + 4(L+\vare)\bigg)\\
\,\,\,\,\,\,& = &\vare diam(b^{-1}(r_1)).
\end{eqnarray*}
and we have completed the proof of Theorem~\ref{GHWarp}.

Note that in $V_{\vare, L}$ was chosen in (\ref{genV0}) and
 (\ref{eqnV3}).

\ProofEnd

We could also consider manifolds in which the function,
\begin{equation} \label{MoreWarpEqn}
\delta(r) :=V_\infty-\frac{Vol_{n-1}(b^{-1}(r))}{r^{p-1}},
\end{equation}
decreases at given rate.  This would give us results
which are stronger than those implied by minimal volume growth
but weaker than those implied by stongly minimal volume growth.

\begin{remark} \label{MoreWarp1} {\em
If we assume that $\delta(r)$ decreases sufficiently
fast, regions of increasing length, like $b^{-1}(r, 2r)$,
 could be shown to be 
Gromov-Hausdorff close to warped product manifolds.  
This can be seen because
we know that for any fixed set of parameters, we can choose
$\delta(r)$ such that Gromov Hausdorff estimating function satisfies
\begin{equation}
\Psi(\delta(r)|n, f, r, 2r, \alpha, \alpha', \xi,0, Dr)
<\vare.
\end{equation}
Note that we are using at most linear diameter growth here 
to say that the diameter of $b^{-1}(r,2r)$ is less than $Dr$,
for some constant, $D$.  
Note also that we do not
bother to rescale the manifold.   Such a theorem would tell us that
\begin{equation}
d_{GH}(b^{-1}(r,2r), b^{-1}(r)\times_f(r,2r)) < \vare.
\end{equation}
In particular, we would truely see the warping of such a manifold.}
\end{remark}

\begin{remark} \label{MoreWarp2}{\em
If we were to take a function, $\delta(r)$, that decreased even
faster, the manifold could be shown to be close to a unique warped product 
manifold.  That is, if we choose $\delta(r)$ such that
for all $r>r_0$ we have
\begin{equation}
\Psi(\delta(r)|n, f, a=r, b=2r, \alpha, 
\alpha', \xi,0, Dr)
<  \vare \left(\frac{1}{2r}\right),
\end{equation}
then 
\begin{equation}
d_{GH}\left(b^{-1}(r,2r), b^{-1}(r)\times_f(r,2r)\right) 
<\vare \left(\frac{1}{2r}\right) .
\end{equation}
Since our Gromov-Hausdorff map is from level set to level set,
if we let $\frac{b^{-1}(r)}{f(r)}$ denote a level set with a localized metric
rescaled by the warping function $f(r)$, then
\begin{equation}
d_{GH}\left(\frac{b^{-1}(r)}{f(r)},\frac{b^{-1}(2r)}{f(2r)}\right) 
< \vare \left(\frac{1}{2r}\right).
\end{equation}
Then, for all $k$,
\begin{equation} 
d_{GH}\left(\frac{b^{-1}(r_0)}{f(r_0)},\frac{b^{-1}(2^Nr_0)}{f(2^Nr_0)}\right)
< \vare \sum_{k=1}^{N+1} \left(\frac{1}{2^kr_0}\right)<\frac{\vare}{r_0}.
\end{equation}
Thus,
\begin{equation}
d_{GH}\left(b^{-1}(r_0,\infty), b^{-1}(r_0)\times_f(r_0,\infty)\right) 
< \vare
\end{equation}
and the manifold is close to a unique warped product manifold.}
\end{remark}

\section {Linear Volume Growth and Sublinear Diameter Growth} \label{Chap7}

\Pind
In this section, we show that a manifold with nonnegative Ricci curvature 
everywhere and linear volume growth has sublinear diameter growth
[Theorems~\ref{SublinDiam} and~\ref{SublinDiam2}].
We prove this both for diameters measured in the ambient manifold and
for localized diameters as defined below in Definition~\ref{DefnDiam}. 

In the previous section, we proved that in such a manifold,
the region $b^{-1}((r, r+L))$ is almost an isometric product after
rescaling by the diameter and taking $r$ large.  Thus the ``diameters'' 
of $b^{-1}(r)$ and  $b^{-1}(r+L)$ are close but only after rescaling
by the diameter of $b^{-1}(r)$.  To get sublinear diameter growth, we need to 
control increasingly long regions, $b^{-1}((r/2, 2r))$ which allows us to
compare $b^{-1}(r/2)$ to $b^{-1}(2r)$.
Here we will rescale by dividing out by $r$  and we use the
fact that $diam(b^{-1}(r)) \le C_D r$ where $C_D$ is the
constant from Theorem~\ref{CompactLevels}, before applying 
Theorem~\ref{ChCo4.85} of [ChCo] and Lemmas
~\ref{MinVolToAlmostMax},~\ref{Rescale},
and~\ref{ThinGH} from Section 2.3.  

Recall the definitions of {\em localized distance} 
and {\em localized component}
from the Section 2.1 [Defn~\ref{DefnDist}, [Note~\ref{components}]] . 

Throughout this chapter, ``diameter'' and $diam$, stated alone,
refer to the diameter measured in the ambient manifold.
We will now define the {\em localized diameter} or, more
precisely, the {\em $s$-almost intrinsic diameter}.

\begin{definition} \label{DefnDiam} {\em
Given any $s>0$, let $U=b^{-1}(R-s, R+s)$.
The {\em $s$-almost intrinsic diameter} or {\em localized
diameter} of a level set of a Busemann function is
\begin{equation}
diam_s(b^{-1}(R)) = \max \{ diam_V(V\cap b^{-1}(R))): 
V \textrm{ is a conn comp of } U \},
\end{equation} 
where }
\begin{equation}
diam_V(V\cap b^{-1}(R))=\sup_{\,\,x,y\in V\cap b^{-1}(R)\,\,} 
                   \inf_{\,c([0,1])\subset V \atop
                                    c(0)=x, c(1)=y} Length(c). 
\end{equation}
\end{definition}

This kind of diameter has been analyzed by Abresch and Gromoll [AbGl].
They proved that the almost intrinsic diameters of
distance spheres in a manifold with nonnegative Ricci curvature grow linearly.
This almost intrinsic diameter of the sphere of radius $R$ was defined to
be the diameter of the largest component of the sphere measured with respect to
a localized distance function $d_U$ where $U=Ann_{R(1-\vare),R(1+\vare)}$.
Thus the almost intrinsic diameter was measured in terms of increasingly thick
annuli.  

In Theorem~\ref{SublinDiam} below, our almost intrinsic diameter 
of Busemann levels is also measured in this way.  However, we prove that
the diameter of the Busemann levels grows sublinearly and we assume
that the manifold has linear volume growth.

\begin{theorem} \label{SublinDiam}
Let $M^n$ be a manifold with nonnegative Ricci curvature and linear
volume growth.
Then given any $\delta \in (0,1)$, 
we have sublinear almost intrinsic diameter growth,
\begin{equation}
\lim_{R\to\infty} \frac {diam_{\delta R}(b^{-1}(R))}{R}=0.
\end{equation}
\end{theorem}

After proving this theorem, we conclude the paper by showing that
manifolds with linear volume growth and nonnegative Ricci curvature
have sublinear diameter growth as well [Theorem~\ref{SublinDiam2}].  
This final theorem does not follow directly from the sublinear 
almost intrinsic diameter 
growth because of the lack of control on the number of localized components 
of the level sets.  

\begin{lemma}  \label{diamgh}
Given any  $\psi>0$ there exists $R_\psi>0$ such that
for all $r \ge R_\psi$, there exists a length space, $X_r$, such that
\begin{equation}
d_{GH}(b^{-1}([r/2, 3r]), X_r \times [r/2,3r])\, \le \,\psi \,r 
\end{equation}
where the localized distance function on $b^{-1}([r/2, 3r])$ is $d_W$
where 
\begin{equation}
W=b^{-1}((r/2-\psi r/2, 3r+\psi r/2)).
\end{equation}

Furthermore the Gromov Hausdorff Equivalence map from
$b^{-1}([r/2, 3r])$ to $X_r \times [r/2,3r]$
is $(\pi(x), b(x))$ where $\pi:b^{-1}(s)\mapsto X_r$ is
a Gromov-Hausdorff equivalence map itself for all $s \in [r/2, 3r]$
with respect to $d_W$.
\end{lemma}

Before proving this lemma, we make a few remarks pointing out some 
useful implications.  Note that we are strongly using the fact that
$X_r \times [r/2,3r]$ is an isometric product. 

\begin{remark} \label{GHMAP1}
Since $F_{GH}$ is almost onto, 
 for all $(x,s)\in X_r \times [r/2,3r]$
there exists $x_s\in b^{-1}([r/2, 3r])\subset M^n$ such that 
$d_{X_r \times [r/2,3r]}(F_{GH}(x_s),(x,s))<\psi r$.
In particular,
$d_{X_r}(\pi(x_s), x)<\psi r$
and $|b(x_s)-s| <\psi r$.  
Furthermore, for
all $y_1, y_2 \in b^{-1}([r/2, 3r])$, 
$d_{X_r}(\pi(y_1), \pi(y_2))< d_W(y_1, y_2) + \psi r$. 
\end{remark}

\begin{remark} \label{smallconn}
In Lemma~\ref{diamgh} we have implicitly stated that the region, 
$W=b^{-1}((r/2-\psi r, 3r+\psi r))$, has the same number of
connected components as $X_r\times (r/2, 3r)$, which is the
same number of componenets as $X_r$ itself.  [Note~\ref{components}]
So if $x$ and $y$ are
in the same connected component of $W$, $\pi(x)$ and $\pi(y)$
are in the same connected component of $X_r$.

On the other hand,  as mentioned in [ChCo],
if $x$ and $y$ are in the same connected component
of $X_r$, then the points $x_s,y_s \in b^{-1}(s-\psi r, s +\psi r)$ 
mentioned in Remark~\ref{GHMAP1}, are in the same
connected component of $b^{-1}(s-4\psi r, s+ 4\psi r)$.  This follows
by dividing up any curve, $C$, between $x$ and $y$ in $X_r$ into
points, $x^i$ of distance less than $\psi r$ apart.  Then there exist
points $x_s^i \in b^{-1}(s-\psi r, s +\psi r)$ as in Remark~\ref{GHMAP1}, 
such that 
\begin{eqnarray*}
d_W(x_s^i, x_s^{i+1}) 
& < & d_{X_r \times [r/2,3r]}
\bigg(F_{GH}(x_s^i),F_{GH}(x_s^{i+1})\bigg)+ \psi r\\
& < & d_{X_r}(\pi(x_s^i), \pi(x_s^{i+1})) + |b(x_s^{i})-b(x_s^{i+1})|+ \psi r\\
& < & d_{X_r}(\pi(x_s^i), x^{i})+d _{X_r}(x^i, x^{i+1})
                 +d_{X_r}(x^{i+1},\pi(x_s^{i+1})) + 3\psi r \\
& < & 3\psi r+3\psi r=6\psi r
\end{eqnarray*}
Thus a piecewise geodesic from $x_s$ to $y_s$ can be drawn
through these points and will remain in the set 
$b^{-1}(s-4\psi r, s+ 4\psi r)$.    
\end{remark}

\noindent {\bf Proof of Lemma~\ref{diamgh} }
Fix any $\psi \in (0,1/2)$.  Let $\omega >0$.  We will
choose the value of $\omega$ later.
Let $R_{\psi, \omega}$ be large enough that
\begin{equation} \label{Rpsiom}
Vol_{n-1}\big(b^{-1}\big(R_{\psi,\omega}(1/2 - \psi)\big)\big) 
> \big(1-\omega\big) V_\infty
\end{equation}
where $V_\infty=\lim_{R\to \infty}Vol_{n-1}(b^{-1}(R))$
as in (\ref{Vinfty}).
Thus, by Lemma~\ref{MinVolToAlmostMax}, 
we have for any $r \ge R_{\psi, \omega}$,
\begin{equation}
U=b^{-1}(r/2 -\psi r, 3r +\psi r)
\end{equation} 
has $\omega$-almost maximal volume [Defn~\ref{AlmostMax}]
with the distance function
\begin{equation}
\rho(x)=d\big(x,b^{-1}(3r +\psi r)\big)
\end{equation} 
and $f(\rho)={\bf{1}}$.
By Lemma~\ref{BuseDist}, $\rho(x)=3r+\psi r -b(x)$.
If we rescale this region by dividing the metric by $r^2$, it
still has $\omega$-almost maximal volume
by Lemma~\ref{Rescale}.

The rescaled region, $U$, can be described as the region between
two level sets of a distance function, $\bar{\rho}(x)$
which is the rescaling of $\rho(x)$, as follows: 
\begin{equation}
U=\bar{\rho}^{-1}(0, 5/2 + 2\psi).
\end{equation} 
By Lemma~\ref{RayDiam}, we have
\begin{eqnarray*}
diam_{rescaled}(U)&=&\frac{diam(b^{-1}(r/2 -\psi r, 3r +\psi r)}{r}\\
&\le& \frac{2(5/2+2\psi)r + diam(b^{-1}(3r+\psi r))}{r}
\le 5 +4 \psi +C_D 
\end{eqnarray*}
where $C_D$ is the diameter growth bound in Note~\ref{LinDiam}.
Now set the parameters $a=0$, $b=(5/2 + 2\psi)$, 
$f(\bar{\rho})={\bf{1}}$, $\Lambda=0$, $D=(5 +4 \psi +C_D) $,
$\alpha=\psi$, $\alpha'=\psi/2$, and $\chi=\psi/4$
and apply the Cheeger-Colding Theorem~\ref{ChCo4.85}.

Thus there exists a function, 
\begin{equation}
\Psi(\omega)
=\Psi(\omega|n,{\bf{1}},a,b,\alpha,\alpha',\chi,\Lambda, D),
\end{equation}
which converges to $0$ as $\omega$ converges to $0$, such that
there exists
a length space $X_{\psi,U}$, depending on
$\alpha'=\psi/2$, and $\chi=\psi/4$ and $U$, such that
\begin{equation}
d_{GH}(\bar{\rho}^{-1}(\psi, 5/2 + \psi), 
X\times (\psi, 5/2 + \psi)) < \Psi(\omega) 
\end{equation}
where the metric on $r^{-1}(\psi, 5/2 + \psi)$ is
$d_W$ where $W=r^{-1}(\psi/2, 5/2 +3 \psi/2)$.
By Remark~\ref{DefnFGH},
the Gromov-Hausdorff map has the form $F_{GH}(x)=(\pi(x),\bar{\rho}(x))$
and by Remark~\ref{GHOnto}, $\pi:\bar{\rho}^{-1}(t)\mapsto X_{\psi,U}$
is also a Gromov-Hausdorff map.  

We now choose $\omega_\psi$ small enough that $\Psi(\omega_\psi) <\psi$.
Rescaling the information back up to full size,
we know there exists $R_\psi= R_{\psi, \omega_\psi}$ of (\ref{Rpsiom})
such that for all $r\ge R_\psi$ there exists $X_r$ such that
\begin{equation}
d_{GH}\big(\rho^{-1}(\psi r, (5/2) r + \psi r),
X_r\times (\psi r, (5/2)\, r + \psi r)\big)        \,  < \,\psi \,r.
\end{equation}
Here the metric is rescaled to 
$d_W$ where $W=\rho^{-1}(\psi/2r, 5r/2 +3r \psi/2 )$.
The Gromov-Hausdorff map has the form $(\pi(x),\rho(x))$
and $\pi:\rho^{-1}(t)\mapsto X_r$ is also a Gromov-Hausdorff map.
This easily implies the lemma using the fact that
$b(x)=3r +\psi r -\rho(x)$.
\ProofEnd


Using the above lemma, we can now compare the localized diameters
of nearby level sets.

\begin{lemma} \label{localdiam}  
Let $M^n$ have nonnegative Ricci curvature and linear volume growth.
Fix $\delta \in (0, 1/2]$.
Given any  $\psi \in (0,\delta/10)$, 
if $r \ge R_\psi$ of Lemma~\ref{diamgh}
then for any $r_1, r_2 \in [r,2r]$ we have, 
\begin{equation}
|diam_{\delta r_1}(b^{-1}(r_1)) - diam_{\delta r_2}(b^{-1}(r_2)| 
\,<\, 6\, \psi \, N_\delta \,r
\end{equation}
where $N_\delta$ is a uniform upper
bound on the number of points $x_i \in X_r$ such that
$d(x_i, x_j) \ge (\delta - 4 \psi)r$.
In particular,
\begin{equation}
N_\delta=  \left(\frac{10\, C_D(1 + \delta)}{\delta}\right) ^n,
\end{equation}    
where $C_D$ is the constant bounding diameter growth of Theorem~\ref{LinDiam}.
\end{lemma}

\noindent{\bf Proof:}                 

By Lemma~\ref{diamgh} we know that
\begin{equation}\label{newgh}
d_{GH}(b^{-1}([r/2, 3r]), X_r \times [r/2, 3r]) < \psi r.
\end{equation}

We first prove that $N_\delta$ is a uniform upper
bound on the number of points $x_i \in X_r$ such that
$d(x_i, x_j) \ge (\delta - 4 \psi)r$.  

Let $x_1, x_2...x_N$ be a maximal set of such points in $X_r$.
By Remark~\ref{GHMAP1}, there exists 
$y_i\in b^{-1}((r-\psi r, r+\psi r))$ such that 
$d_{X_r}(\pi(y_i), x_{i})<\psi r$.  In particular,
\begin{eqnarray}
d_W(y_i, y_j) &>& d_{X_r\times [r/2, 3r]}(F_{GH}(y_i),F_{GH}(y_j)) -\psi r\\
& \ge & d_{X_r}(\pi(y_i), \pi(y_j)) - \psi r\\
& \ge & d _{X_r}(x_i, x_j)
      - d_{X_r}(\pi(y_i), x_{i})-d_{X_r}(x_j,\pi(y_j)) -\psi r \\
& > & (\delta - 4 \psi)r -3\psi r \,\, = \,\, (3 \delta/10)r.
\end{eqnarray}
Thus $B_{y_i}(\delta r/10)$ are disjoint balls in 
$b^{-1}(r-\delta r/5, r+\delta r/5)$.  

Recall the linear diameter growth constant, $C_D \ge 1$,
 of Theorem~\ref{LinDiam}.
Since $y_i$ are in $b^{-1}(r-\delta r/5, r + \delta r/5)$, which
has diameter less than or equal to $2(\delta r/5) + C_D(r + \delta r/5)$
by Lemma~\ref{RayDiam}, we can apply the
Relative Volume Comparison Theorem to bound the
volumes of the balls from below.
\begin{eqnarray*}
Vol\big(B_{y_i}(\delta r/10)\big) & \ge &
\frac {\big(\delta r/10\big)^n 
           Vol\big(B_{y_i}(C_D(r+\delta r/5) +2\delta r /5)\big)}
       {\big(C_D(r+\delta r/5)+2\delta r /5\big)^n} \\
& \ge & \left(\frac {( \delta r/10)}{C_D(r+\delta r)}\right)^n
Vol(b^{-1}(r-\delta r/5, r + \delta r/5)).
\end{eqnarray*}
Since the balls are all disjoint,
\begin{eqnarray*}
Vol\big(b^{-1}((r-\delta r/5, r + \delta r/5))\big)
&\ge & \sum_{i=1}^N \,\,Vol\big(B_{y_i}(\delta r/5)\big)\\
&\ge& N  \frac {(\delta)^n Vol(b^{-1}(r-\delta r/5, r + \delta r/5))}
                {(10\,C_D(1+\delta))^n},
\end{eqnarray*}
and we have the uniform  estimate for $N_\delta$.


We now prove that for all $s \in [r, 2r]$,
\begin{equation} \label{diampart1}
diam_{\delta s}(b^{-1}(s)) < diam(X_r) + 3 \psi N_\delta r,
\end{equation}
where $diam(X_r)$ is the diameter of the largest connected
component of $X_r$.  

Recall Definition~\ref{DefnDiam} of
the almost intrinsic diameter.  Let $U=b^{-1}(s-\delta s, s-\delta s)$
and let $V$ be the largest pathwise connected component of $U$;
so 
\begin{equation}\label{diam1n2}
diam_{V}(b^{-1}(s) \cap V)=diam_{\delta s}(b^{-1}(s)).
\end{equation}
Let $x$ and $y$ be any pair of points in $b^{-1}(s) \cap V$.
We claim that for all $\vare >0$, there exists a curve, $c$, 
contained in $V$ from $x$ to $y$ such that 
\begin{equation}\label{claim1}
L(c) <diam(X_r) + 3N_\delta \psi r +\vare.
\end{equation}
Once we have proven the claim, (\ref{diampart1}) is proven.

Fix $\vare$.  To find a curve, $c$, we first note that 
$x$ and $y$ are in the same connected component of $b^{-1}(r/2,3r)$.
Thus $\pi(x)$ and $\pi(y)$ are in the same connected componenet of
$X_r$.  So there exists a curve, $C_\vare$,  in the length space
$X_r$ between $\pi(x)$ and $\pi(y)$ such that 
\begin{equation}
d_{X_r}(\pi(x),\pi(y))\le L(C_X)\le d_{X_r}(\pi(x),\pi(y))+\vare.
\end{equation}
Let $p_0=\pi(x), p_1, p_2,...p_N=\pi(y)$ be equally spaced
points along $C_\vare$ such that $N = N_\delta+1$.  Thus 
$d(p_i, p_{i+1})<(\delta - 4 \psi)r$ and 
\begin{equation}
\sum_{i=0}^{N-1} d_{X_r}(p_i, p_{i+1}) 
\le L(C_\vare).
\end{equation} 

By Remark~\ref{GHMAP1}, there exists 
$x_i\in b^{-1}((s-\psi r, s+\psi r))$ such that
$x_0=x$, $x_N=y$ and $d_{X_r}(\pi(x_i), p_i)< \psi r$.  Thus,
\begin{eqnarray*}
d(x_i, x_{i+1})& \le & d_W(x_i, x_{i+1})
\,\,< \,\,d_{X_r}(\pi(x_i), \pi(x_{i+1})) +\psi r \\
&< & 3 \psi r + d_{X_r}(p_i, p_{i+1}) \,\,<\,\,\delta r -\psi r
\,\,<\,\,\delta s,
\end{eqnarray*}
where $W=b^{-1}(r/2 -\psi r/2, 3r +\psi r/2)$ as in Lemma~\ref{diamgh}. 
So a minimal geodesic from $x_i$ to $x_{i+1}$ is contained in 
$b^{-1}((s-\delta s, s+ \delta s))\subset W$.  So there is a piecewise minimal
geodesic curve from $x$ to $y$ contained in $b^{-1}((s-\delta s, s+ \delta s))$
of length,
\begin{eqnarray*}
\sum_{i=0}^{N-1} d_W(x_i, x_{i+1}) &\le& 
\sum_{i=0}^{N-1} d_{X_r}(p_i, p_{i+1}) + N_\delta 3\psi r \\
&\le&  L(C_\vare) + 3 N_\delta \psi r \\
& \le & diam(X_r) + \vare + 3 N_\delta \psi r. 
\end{eqnarray*}
This piecewise geodesic satisfies our claim, (\ref{claim1}), 
and (\ref{diampart1}) follows.


We now prove that for all $s \in [r, 2r]$,
\begin{equation}\label{diampart2}
diam(X_r) < diam_{\delta s}(b^{-1}(s))  + 3\psi  r.
\end{equation}

First recall that $diam(X_r)$ is the diameter of its largest connected
component, $Z_r\subset X_r$.  Let $x$ and $y$ be any pair of points in $Z_r$.
There exists $p, q \in b^{-1}((s-\psi r, s+\psi r))$ such that
\begin{equation} \label{toodle}
d(\pi(p),x)<\psi r \textrm{  and  } d(\pi(q),y)<\psi r,
\end{equation}
as mentioned in Remark~\ref{GHMAP1}.

Now, $p$ and $q$ are in the same connected component
of $b^{-1}((s-\delta s, s +\delta s))$ by Remark~\ref{smallconn} and 
by the fact that $b^{-1}((s-\delta s, s +\delta s))$ contains the region 
$b^{-1}((s-4 \psi r, s+ 4 \psi r))$. 
So for all $\vare >0$, there exists a curve 
$C_\vare\in b^{-1}((s-\delta s, s+\delta s))\subset W$ which almost achieves
the $s\delta$-almost intrinsic distance between $p$ and $q$ and
\begin{equation}\label{oodle1}
d_W(p,q) \le L(C_\vare) \le d_{s\delta}(p,q) +\vare 
\le diam_{\delta s}(b^{-1}(s)) +\vare.
\end{equation}
On the other hand, by (\ref{newgh}),(\ref{toodle}) and 
the triangle inequality,
\begin{eqnarray}\label{oodle2}
d_W(p,q)&\ge &d_{X_r\times[r/2,3r]}(F_{GH}(p),F_{GH}(q)) -\psi r\\
&\ge& d_{X_r}(\pi(p), \pi(q))-\psi r \ge d_{X_r}(x,y)- 3\psi r.
\nonumber
\end{eqnarray}
Combining (\ref{oodle1}) and (\ref{oodle2}), taking $\vare \to 0$ and then 
maximizing over all $x$ and $y$ in $Z_r$,
we have (\ref{diampart2}).

It is easy to see that (\ref{diampart1}) and (\ref{diampart2}) applied
to $s=r_1$ and $s=r_2$ alternatively, imply the lemma.
\ProofEnd
 
We can now prove the sublinear almost
intrinsic diameter growth theorem.

\noindent {\bf Proof of Theorem~\ref{SublinDiam}:}
We need to show that given any $\delta \in(0,1/2)$,  
we have sublinear almost intrinsic diameter growth,
\begin{equation}
\lim_{R\to\infty} \frac {diam_{\delta R}(b^{-1}(R))}{R}=0.
\end{equation}

Given any $\vare >0$, let 
\begin{equation}
\psi=\min \{\delta/10, \vare/3 N_\delta \}.
\end{equation} 
Thus by Lemmas~\ref{diamgh} and~\ref{localdiam},
there exists $R_\psi$ such that for all $r \ge R_\psi$
and for all $r'\in [r,2r]$,
\begin{equation}\label{rto2rhere}
diam_{\delta r'}(b^{-1}(r')) < \vare r + diam_{\delta r}(b^{-1}(r)).
\end{equation}
Thus, applying this repeatedly to
$r=2^j \delta R_\vare$ and $r'=2r$, we have
\begin{eqnarray*}
diam_{2^k \delta R_\vare}(b^{-1}(2^k R_\vare) 
& < & \vare 2^{k-1}R_\vare 
             + diam_{\delta 2^{k-1}R_\vare}(b^{-1}(2^{k-1}R_\vare))\\
& < & \vare 2^{k-1}R_\vare + \vare 2^{k-2}R_\vare 
             + diam_{\delta 2^{k-2}R_\vare}(b^{-1}(2^{k-2}R_\vare))\\
&< & \vare (2^{k-1} + ... +2 +1) R_\vare
           + diam_{\delta R_\vare}(b^{-1}(R_\vare))\\
& < & \vare 2^k R_\vare + diam_{\delta R_\vare}(b^{-1}(R_\vare))
\end{eqnarray*}
For all $R> R_\vare$ there exists $k$ such that 
$R \in [2^kR_\vare, 2^{k+1}R_\vare]$ so by (\ref{rto2rhere}), we have
\begin{equation}
\frac{\,diam_{\delta R}(b^{-1}(R))\,}{R} \,\,< \,\,
\vare  \,+ \,
2 \, \frac{\,diam_{2^k \delta R_\vare}(b^{-1}(2^k R_\vare))\,}{2^kR_\vare} 
\end{equation}
Combining this with the above estimate and taking $R$ to infinity,
we get
$$ 
\limsup_{R\to\infty}  \frac{diam_{\delta R}(b^{-1}(R))}{R} \le 
\vare + 2 \lim_{k\to\infty}
\frac{\,\vare 2^k R_\vare + diam_{\delta R_\vare}(b^{-1}(R_\vare))\,}
{2^k R_\vare} = 3 \vare.
$$
Since this is true for all $\vare >0$ we have sublinear diameter growth.
\ProofEnd

%
%

We now prove our final theorem, Theorem~\ref{SublinDiam2}, that
\begin{equation}
\lim_{R\to\infty} \frac {diam(b^{-1}(R))}{R}=0.
\end{equation}

{\bf Proof of Theorem~\ref{SublinDiam2}:}
By Theorem~\ref{LinDiam} of [So2], we already know that
\begin{equation}
\frac {diam(b^{-1}(R))}{R}\le C_D < \infty.
\end{equation}

If we assume that the diameter growth is not sublinear, then
there exists a sequence, $r_i$, approaching infinity such that
\begin{equation} \label{CLdiam}
\frac {diam(b^{-1}(r_i))}{r_i} \ge C_L>0 \qquad \forall i.
\end{equation} 
So there exist $x_i, y_i \in b^{-1}(r_i)$, and there exists, $\sigma_i$,
a minimal geodesic from $x_i$ to $y_i$ such that $L(\sigma_i)=h_ir_i$
where $h_i \in [C_L,C_D]$.

Suppose there is a subsequence, $i_j$, such that 
\begin{equation} \label{diamnearby}
\sigma_{i_j} \subset b^{-1}([r_{i_j}/2, 3r_{i_j}/2]).
\end{equation}
Then we have a minimal geodesic of length $h_{i_j}r_{i_j}$
contained in this region, so 
\begin{equation}
diam_{\frac 1 2 (r_{i_j})}(b^{-1}(r_{i_j})) \ge h_{i_j}r_{i_j}\ge C_Lr_{i_j}.
\end{equation}
This contradicts the sublinear 1/2-almost intrinsic diameter growth of
the manifold [Theorem~\ref{SublinDiam}].
Thus there exists $N$ such that for all $i \ge N$,
$\sigma_{i}$ is not a subset of $b^{-1}(r_{i}/2, 3r_{i}/2)$.

Let 
\begin{equation}\label{SigmaMin}
S_i\,=\,\min_{t\in [0,h_ir_i]} \,b(\sigma_i(t)).
\end{equation} 
and let
\begin{equation}\label{TMax}
T_i=\max_{t\in [0,h_ir_i]} b(\sigma_i(t)).
\end{equation}
For all $i \ge N$ either $2 S_i <r_{i}$ or $2T_i/3 >r_{i}$. 

Suppose there is a subsequence, $i_j$, such that $2T_{i_j}/3 >r_{i_j}$.
Then $\sigma_{i_j}$ is a minimal geodesic which starts at $x_{i_j}$ in
$b^{-1}(r_{i_j})$, passes through a point, $x'_{i_j}$ in $b^{-1}(2T_{i_j}/3)$,
passes through a point in $b^{-1}(T_{i_j})$, continues back through
a point, $y'_{i_j}$ in $b^{-1}(2T_{i_j}/3)$, before returning to $y_{i_j}$ in
$b^{-1}(r_{i_j})$.  Thus there are points $x'_{i_j}$ and $y'_{i_j}$ in
$b^{-1}(2T_{i_j}/3)$ with a minimal geodesic between them of length at least
$2T_{i_j}/3$ which remains in $b^{-1}(2T_{i_j}/3, T_{i_j})$.  So
\begin{equation}\label{contradiam1}
diam_{\frac 1 2 (2T_{i_j}/3)}(b^{-1}(2T_{i_j}/3)) \ge 2T_{i_j}/3.
\end{equation}
As $r_{i_j}$ goes to infinity, $T_{i_j}$ approaches infinity,
and then (\ref{contradiam1}) contradicts Theorem~\ref{SublinDiam} with
$\delta =1/2$.
Thus there exists $N'$ such that $2T_{i}/3 \le r_{i}$ for all $i \ge N'$;
so $2 S_i <r_{i}$ for all $i \ge N'$.

We would like to use the same trick with the $S_i$
of (\ref{SigmaMin}). but first we
must show that $S_i$ diverge to infinity.

Suppose there exists $R$, and there exists a subsequence $i_j$ such that
$S_{i_j} \le R$.  Thus there exists $t_{i_j}$ such that 
$\sigma_{i_j}(t_{i_j})\in b^{-1}(R)$.  Since this level set is compact,
a subsequence of the $\sigma_{i_j}$ must converge.  Since
$L(\sigma_{i_j}) \ge C_L r_{i_j}$, this subsequence must converge to
a line.  So by the Splitting Theorem of Cheeger and Gromoll, the
manifold is split [ChGl].  However, this implies that $b^{-1}(r_i)$ is totally 
geodesic and so $S_i=r_i$, which contradicts $2S_i < r_i$.

Thus $S_i$ diverges to infinity and $S_i \le r_i/2$.
Then $\sigma_{i}$ is a minimal geodesic which starts at $x_{i}$ in
$b^{-1}(r_{i})$, passes through a point, $x'_{i}$ in $b^{-1}(2S_{i})$,
passes through a point in $b^{-1}(S_{i})$, continues back through
a point, $y'_{i}$ in $b^{-1}(2S_{i})$, before returning to $y_{i}$ in
$b^{-1}(r_{i})$.  Thus there are points $x'_{i}$ and $y'_{i}$ in
$b^{-1}(2S_{i})$ with a minimal geodesic between them of length at least
$2S_i$ which remains in $b^{-1}(S_{i}, 2S_{i})$.  So
\begin{equation}
diam_{\frac 1 2 (2S_{i})}(b^{-1}(2S_{i})) \ge 2S_{i},
\end{equation}
which contradicts Theorem~\ref{LinDiam} for $\delta =1/2$ as
$S_{i}$ approaches infinity.

Thus $C_L$ of (\ref{CLdiam}) cannot exist, and the manifold has
sublinear diameter growth.
\ProofEnd

Similar theorems may be provable for manifolds with a quadratically 
decaying lower bound on Ricci curvature, but one must be careful to
rescale $X_r$ as it is compared to each level set.

\end{document}